\documentclass[a4paper, notitlepage, 12pt]{article}

\usepackage[total={160mm,250mm}, top=25mm, left=25mm, includefoot]{geometry}

\usepackage{ucs}              
\usepackage[utf8x]{inputenc} 	

\usepackage[T1]{fontenc}     
\usepackage{lmodern}         

\usepackage[singlespacing]{setspace}

\usepackage{color}
\usepackage{amsmath,amsthm,amssymb,amscd,amsbsy}
\usepackage{upgreek}

\usepackage{tensor}

\usepackage{wasysym}
\usepackage{tocbibind}          
\usepackage{float}
\usepackage{tabularx,paralist}
\usepackage{graphicx}
\usepackage{hyperref}
\usepackage{authblk}

\usepackage{multirow} 
\usepackage{array}
\usepackage{caption}

\usepackage{rcs} \RCS $Revision: 2.7 $ \RCS $Date: 2018/06/16 00:23:47 $ \RCS $Author: hgm $ \RCS $RCSfile: 17_Param-Lin_Map.tex,v $ \RCS $Id: 17_Param-Lin_Map.tex,v 2.7 2018/06/16 00:23:47 hgm Exp $

\usepackage[british]{babel}

\usepackage{refdef}
\usepackage{ifontdef}

\newtheorem{thm}{Theorem}

\newtheorem{prop}[thm]{Proposition}

\newtheorem{coro}[thm]{Corollary}

\newtheorem{defi}[thm]{Definition}

\newtheorem{lem}[thm]{Lemma}

\newcommand{\ignore}[1]{}

\newcommand{\Lp}{\mrm{L}}

\newcommand{\vsigma}{\varsigma}
\newcommand{\vtheta}{\vartheta}
\newcommand{\vphi}{\varphi}
\newcommand{\vpi}{\varpi}
\newcommand{\vkappa}{\varkappa}

\newcommand{\vek}[1]{\mathchoice{\displaystyle\boldsymbol{#1}}
{\textstyle\boldsymbol{#1}}{\scriptstyle\boldsymbol{#1}}
{\scriptscriptstyle\boldsymbol{#1}}}

\newcommand{\mat}[1]{\mathchoice{\displaystyle\mathbf{#1}}
{\textstyle\mathbf{#1}}{\scriptstyle\mathbf{#1}}
{\scriptscriptstyle\mathbf{#1}}}

\newcommand{\opb}[1]{\vek{{\mathsf{#1}}}}

\newcommand{\ops}[1]{\mathchoice{\displaystyle\mathsf{#1}}
{\textstyle\mathsf{#1}}{\scriptstyle\mathsf{#1}}
{\scriptscriptstyle\mathsf{#1}}}

\newcommand{\tnb}[1]{\mathchoice{\displaystyle\mathboldsans{#1}}
{\textstyle\mathboldsans{#1}}{\scriptstyle\mathboldsans{#1}}
{\scriptscriptstyle\mathboldsans{#1}}}

\newcommand{\tns}[1]{\mathchoice{\displaystyle\mathsans{#1}}
{\textstyle\mathsans{#1}}{\scriptstyle\mathsans{#1}}
{\scriptscriptstyle\mathsans{#1}}}

\newcommand{\EXP}[1]{\mathbb{E}\left(#1\right)}

\newcommand{\diag}{\mathop{\mathrm{diag}}\nolimits}

\newcommand{\im}{\mathop{\mathrm{im}}\nolimits}
\newcommand{\dom}{\mathop{\mathrm{dom}}\nolimits}

\newcommand{\tr}{\mathop{\mathrm{tr}}\nolimits}

\newcommand{\spn}{\mathop{\mathrm{span}}\nolimits}

\newcommand{\di}{\mathrm{d}}

\newcommand{\ii}{\mathchoice{\displaystyle\mathrm i}
{\textstyle\mathrm i}{\scriptstyle\mathrm i}
{\scriptscriptstyle\mathrm i}}

\newcommand{\KL}{Karhunen-Lo\`eve}
\newcommand{\ip}[2]{\langle #1, #2 \rangle}
\newcommand{\bkt}[2]{\langle #1 | #2 \rangle}

\newcommand{\nd}[1]{\| #1 \|}

\definecolor{myred}{rgb}{1, 0.2, 0.2}

\newcommand{\autheadcr}{\authorcr}
\newcommand{\citep}[1]{\cite{#1}}

\newcommand{\authorhgm}{Hermann G. Matthies}
\newcommand{\authorro}{Roger Ohayon}

\newcommand{\affilwire}{Institute of Scientific Computing\autheadcr
                        Technische Universit\"at Braunschweig\autheadcr
                        38092 Braunschweig, Germany\autheadcr
                        e-mail: \ttt{wire@tu-bs.de}}
\newcommand{\affilcnam}{Laboratoire de M\'ecanique des Structures et des
                       Syst\`emes Coupl\'es\autheadcr
                       Conservatoire National des Arts et M\'etiers (CNAM)\autheadcr
                      75141 Paris Cedex 03, France}

\newcommand{\thetitle}{Analysis of parametric models\\
                        --- linear methods and approximations}

\newcommand{\theauthor}{\authorhgm, \authorro}
\newcommand{\thesubject}{35B30, 37M99, 41A05, 41A45, 41A63, 60G20, 60G60, 65J99, 93A30}
\newcommand{\thekeywords}{parametric models, reproducing kernel Hilbert space,
                          correlation, factorisation, spectral decomposition, representation}

\newcommand{\textdate}{}

\newcommand{\thecommon}{.}
\newcommand{\thebib}{./bib}

\include{config}

\include{\thecommon/default-commands}

\begin{document}

\title{\thetitle\thanks{Partly supported by the Deutsche
          Forschungsgemeinschaft (DFG) through SFB 880 and SPP 1886.}}


\author[*]{\authorhgm}
\author[o]{\authorro}


\affil[*]{\affilwire}
\affil[o]{\affilcnam}

\date{\textdate}


\ignore{          


\setcounter{page}{0}
\thispagestyle{empty}
\cleardoublepage

\include{titlepage}

\newpage

\thispagestyle{empty}
\vspace*{\stretch{2}}

\begin{flushleft}
\begin{tabular}{ll}
\makeatletter
This document was created \textdate{} using \LaTeXe. \\[1cm]
\makeatother
\end{tabular}

\begin{tabular}{ll}
\begin{minipage}{6cm}
Institute of Scientific Computing\\ 
Technische Universit\"at Braunschweig\\
M\"uhlenpfordstra\ss{}e 24\\
D-38106 Braunschweig, Germany\\

\texttt{url: \url{www.wire.tu-bs.de}}\\
\makeatletter
\texttt{mail: \href{mailto:wire@tu-bs.de?subject=\thetitle}{wire@tu-bs.de}}
\makeatother
\end{minipage}
&
\begin{minipage}{2.5cm}
\vspace{-0.5cm}
\includegraphics[width=2.4cm]{common/logo_wire_ohnekreis}

\end{minipage}
\end{tabular}

\vspace*{\stretch{1}}

Copyright \copyright{} by \theauthor{}\\[5mm]
\end{flushleft}

This work is subject to copyright. All rights are reserved, whether the whole or part of the material is concerned, specifically the rights of translation, reprinting, reuse of illustrations, recitation, broadcasting, reproduction on microfilm or in any other way, and storage in data banks. Duplication of this publication or parts thereof is permitted in connection with reviews or scholarly analysis. Permission for use must always be obtained from the copyright holder.\\[5mm]

Alle Rechte vorbehalten, auch das des auszugsweisen Nachdrucks, der auszugsweisen oder vollständigen Wiedergabe (Photographie, Mikroskopie), der Speicherung in Datenverarbeitungsanlagen und das der Übersetzung.


}            

\maketitle

%

\begin{abstract}
Parametric models in vector spaces are shown to possess an associated linear map.
This linear operator leads directly to reproducing kernel Hilbert spaces and
affine- / linear- representations in terms of tensor products.  From the associated
linear map analogues of covariance or rather correlation operators can be formed.
The associated linear map in fact provides a factorisation of the correlation.
Its spectral decomposition, and the associated \KL{}- or proper orthogonal decomposition
in a tensor product follow directly.   It is shown that all factorisations of a
certain class are unitarily equivalent, as well as that every factorisation
induces a different representation, and vice versa.

A completely equivalent spectral and factorisation analysis can be carried out in
kernel space.  The relevance of these abstract constructions is shown on a number
of mostly familiar examples, thus unifying many such constructions under
one theoretical umbrella.  From the factorisation one obtains tensor representations,
which may be cascaded, leading to tensors of higher 
degree.  When carried over to a discretised level in the form of a model order reduction,
such factorisations allow very sparse low-rank approximations which lead to very
efficient computations especially in high dimensions.

\vspace{5mm}

\vspace{5mm}
{\noindent\textbf{MSC subject classifications:} \thesubject}


\end{abstract}

%
%
%
%
%
%
%
%
%
%
%



\cleardoublepage

\tableofcontents
\cleardoublepage





%

\section{Introduction}  \label{S:intro}
Parametric models are used in many areas of science, engineering, and
economics, etc.  They appear in cases of \emph{design} of some systems,
where the parameters may be design variables of some kind, and the variation
of the parameters may show different possibilities, or display the
envelope of the system over a range of possible parameter values.  Other
possibilities arise when one wants to \emph{control} the behaviour of
some system, and the parameters are control variables.  This is
closely connected with situations where one wants to \emph{optimise}
the behaviour of some system more generally by changing some parameters.
Another important area is where some of the parameters may be uncertain,
i.e.\ they could be random variables, and with respect to these one
wants to perform \emph{uncertainty quantification}.  Of course it is also
possible that the parameter set has many-fold purposes, for example that some of
the parameters model design variables, while others are uncertain,
cf.\ \citep{SoizeFarhat2017}.

Often the problem of having to deal with a parametric system is compounded
by the fact that one also has to approximate the system behaviour through
a \emph{reduced order model} due to high computational demands of the
full system.  This reduced model then therefore becomes a parametrised
reduced order model.  The survey \citep{BennWilcox-paramROM2015} and
the recent collection \citep{MoRePaS2015} as well as
the references therein provide a good view not only of reduced
order models which depend on parameters, but also of parametric problems
in general and some of the areas where they appear.  So for further
information on parametrised reduced order models and how to generate
them we refer to these references.

Here, we want to concentrate on a
certain topic illuminating the theoretical background of such parametrised models.
This is the connection between separated series representations, associated
linear maps, the singular value (SVD) and proper orthogonal decomposition (POD),
and tensor products.  This then immediately opens the connections to
reduced order models and low-rank tensor approximations.  It will be seen that
the distribution of singular values of the associated linear map determines
how many terms are necessary for a good approximation with a reduced order model.
For higher order tensor representations in the context of hierarchical tensor
approximations, it is the SVD structure of the tensor product splittings
associated with the tree structure of the index set partitions.

Typically, the parameters are assumed to be tuples of
independent real numbers, but here no assumptions are made about the parameter set.
The geometry and topology of the parameter set is reflected by set of
real functions defined on the parameter set, which can be viewed like co-ordinates
in this context.

In some cases, like design evaluations and uncertainty quantification, the
parameter set in itself is not important, but only the range or distribution
of the parametric model.  Here the analysis of the associated linear map
allows the re-parametrisation of the parametric model with the parameters
taken from a different set.  The principal result is then that within a certain
class of representations of parametric models and associated linear maps there
is a one-to-one correspondence between separated series representations
and factorisations of the associated linear map.  But in general the representation
through a linear map is much more general, and allows the modelling of
\emph{weak} or \emph{generalised} models.

As a possible starting point to introduce the subject, assume
that some physical system is investigated, which
is modelled by an evolution equation for its state: 
\begin{equation} \label{eq:I}
   \frac{\di}{\di t}u(t) = A(q;u(t)) + f(q;t);\quad u(0) = u_0,
\end{equation}
where $u(t) \in \C{U}$ describes the state of the system
at time $t \in [0,T]$ lying in a Hilbert space $\C{U}$ (for the sake of
simplicity), $A$
is an operator modelling the physics of the system, and $f$ is some external
influence (action / excitation / loading).  Assume that the model depends on some
quantity $q$, and assume additionally that for all $q$ of interest the system
\refeq{eq:I} is well-posed.   One part of these parameters $q$ may describe
the specific nature of the system \refeq{eq:I}, whereas another part of
the parameters, here denoted by $p\in\C{P}$ has to be varied
for one reason or another in the analysis.

One is often interested in how the system changes when these
\emph{parameters} $p$ change.  The parameter $p$ can be for example
\begin{itemize}
\item just the quantity, $p=q$; or
\item the quantity and the action, $p=(q,f)$; or
\item as before, but including the initial condition, $p=(q,f,u_0)$; or
\item many other combinations.
\end{itemize}
To deal with all these different possibilities under one notation,
the state \refeq{eq:I} can be rewritten as 
\begin{equation} \label{eq:I-p}
   \frac{\di}{\di t}u(t) = A(p;u(p;t)) + f(p;t);\quad u(0) = u_0,
\end{equation}
with a solution $u(p;t)$ denoting the explicit dependence on the parameter
$p\in\C{P}$.

Frequently, the interest is in functionals of the state $\Psi(p,t)=\Psi(p,u(p;t))$,
and the functional dependence of $u$ on $p$ becomes important.
Such situations arise in design, where $p$ may be a design parameter
still to be chosen, and one may seek a design such that a functional
$\Psi(p,t)$ or some kind of temporal integral or average
$\tns{\psi}(p) = \int_0^T \Psi(p,t) \rho(t)\,\di t$ is,
e.g., maximised \citep{Luenberger}.  Optimal control is a special case of this, 
as one may try to influence the time evolution in such a way that $\Psi(p,T)$
(or $\tns{\psi}(p)$ above)
is minimised or maximised.  Another example is when the $p\in\C{P}$ are uncertain
parameters, modelled by random variables.  In the process of uncertainty
quantification \citep{Matthies_encicl, xiu2010numerical, knio2010spectral, Smith2014}
one then may want to compute expected values $\D{E}_p(\Psi(p,t))$.  It may
also be that the parameters have to be determined or identified to allow
the model to match some observed behaviour, this is called an inverse
problem, see \citep{hgm17-ECM} and the references therein.  Another case
is a general design evaluation, where one is interested in the range of
$u(p;t)$ --- or $\Psi(p,t)$ or $\tns{\psi}(p)$ --- as $p$ varies over $\C{P}$.

The situation just sketched involves a number of objects which are functions
of the parameter values $p$.  While evaluating $A(p)$ of $f(p)$ in \refeq{eq:I-p}
for a certain $p$ may be straightforward, one may easily envisage situations
where evaluating $u(p;t)$ or $\Psi(p,t)$ may be very costly as it may involve
some very time consuming simulation or computation.
Therefore one is interested in representations of $u(p;t)$ or $\Psi(p,t)$
which allow the evaluation in a cheaper way.  These simpler representations
are called by many names, some are \emph{proxy}- or \emph{surrogate}-model.
As will be shown in the following \refS{parametric}, any such
parametric object may be represented in many different ways, many of which
can be analysed by linear maps which are associated with
such representations.  It will be shown that these representations may be
seen as an element of a tensor product space.  This in turn can be used to
find very sparse approximations to those objects, and in turn much cheaper
ways to evaluate the functional $\Psi$ or $\tns{\psi}$ for other parameter values.

This association of parametric models and linear mappings has probably been
known for a long time.  The first steps which the authors are aware of in
this direction of using linear methods are \citep{Loeve1945, Loeve1946}.  A seminal
work  in this kind of inquiry is \citep{Karhunen1947} (with English
\citep{Karhunen1947-e} and Spanish \citep{Karhunen1947-s} translations
available online), which contains a first rather thorough exploration
in the context of probability on infinite dimensional vector or function spaces,
and many influential ideas, see also \citep{Karhunen1946, Karhunen1950}
for similar work.   The name \KL{} expansion for approximations of this
kind in the context of probability theory was coined after the authors
of \citep{Karhunen1947} and \citep{Loeve1945, Loeve1946}.  This name
is used in this context, in other areas that representation is now often known
under the name \emph{proper orthogonal decomposition} (POD), which is
firmly connected with the \emph{singular value decomposition} (SVD)
of the associated linear map.  In following publications,  see
\citep{segal56-TAMS, segal58-TAMS, LGross1962, segalNonlin1969}, the terms
\emph{generalised processes}, \emph{linear process}, \emph{weak distribution}
or equivalently \emph{weak measure} or \emph{weak process} appear
for the associated linear map.  This is indicative of the generalisation
possible with these linear methods, see also the monographs
\citep{gelfand64-vol4, LoeveII, kreeSoize86}.

A first step of reviving and also connecting these methods of analysis
with the theory of low-rank tensor approximations \citep{Hackbusch_tensor}
was undertaken in \citep{boulder:2011} in the context of uncertainty
quantification and inverse problems.  It is furthermore also connected
with non-orthogonal decompositions which are easier to compute,
like the \emph{proper generalised decomposition}(PGD) \citep{chinestaBook}.
Here we continue that endeavour of showing the connection of parametric
models, model reduction of parametric models, and sparse numerical
approximations
to a certain extent in a more general setting.  It is on this theoretical
background that one may analyse modern numerical methods which
allow the numerical handling of very high-dimensional objects,
i.e.\ where one has to deal with an essentially high-dimensional
space for the parameters $p\in\C{P}$.

Whereas the parametric map may be quite complicated, the association
with a linear map translates the whole problem into one of linear
functional analysis, and into linear algebra upon approximation and
actual numerical computation.  Also, whereas
the set $\C{P}$ might have a quite arbitrary structure,
this is replaced by a subspace of the vector space $\D{R}^{\C{P}}$ of real
valued functions on $\C{P}$, in some sense this is a `problem oriented
co-ordinate system' on the set.
This is a frequent technique in mathematics,
and it replaces the quite arbitrary set by a vector space, which is
much more accessible.  Let us recall a situation which is similar
and may be well-known to many readers.  When the need arose to deal
with very singular functions, especially when Dirac needed an `ideal'
object like the $\updelta$-`function', for this and other so-called
\emph{generalised} functions or \emph{distributions} a fruitful 
mathematical formulation turned out to be the model
of a linear map into real numbers on a space of smooth regular
functions, see e.g.\ \citep{gelfand64-vol1, gelfand64-vol2}.

The association with a linear map is quickly shown to be related to
representations connected with the adjoint of the map, and the
precise definition and properties of the associated linear map
are given in \refS{parametric}.  The connection with reproducing
kernel Hilbert spaces (RKHS) \citep{berlinet} also appears naturally
here, and it is shortly sketched.  From the map
and its adjoint we obtain the `correlation', which will be analysed
in \refS{correlat}.  Here the spectral analysis and factorisation
of the correlation will become important 
\citep{Segal1978, reedSimon-vol1, reedSimon-vol2,DautrayLions3}.
This also connects the whole idea of linear methods for representation
with tensor representations, which appear naturally in the spectral
analysis.  The kernel, which on the RKHS is the reproducing kernel, now
appears in another context than the one already alluded to in
\refS{parametric}, and in \refS{kernel} the kernel side of the
representation is analysed, which is the classical domain of integral transforms
and integral equations as already envisaged in \citep{Karhunen1947}.
Some examples and interpretations are explained in \refS{xmpls}, to give an idea of the
breadth of possible applications of the theory.  Here  the connection
of these linear methods to both linear model reduction and nonlinear
model reduction in the form of low-rank tensor approximations
\citep{Hackbusch_tensor} is mentioned and briefly explained.
The last \refS{refine} before the conclusion in \refS{concl} deals
very shortly with certain refinements which are possible when some
a priori structure of the represented spaces is known; we have
connected it here with vector- and tensor-fields.

%
%

%
%
%
%
%
%
%
%
%
%
%
%
%
%
%
%
%



\section{Parametric problems} \label{S:parametric}
Let $r: \C{P} \rightarrow \C{U}$ be a parametric description of
one of the objects alluded to in the introduction, or
the state or
response of some system, where $\C{P}$ is some set, and $\C{U}$ 
is assumed --- for the sake of simplicity ---
as a separable Hilbert space with inner product
$\langle\cdot|\cdot\rangle_{\C{U}}$.  More general locally convex
vector spaces are possible, but the separable Hilbert space is in many
ways the simplest model.

The situation in its purest form may
be thought of in an abstract way as follows:
$ F: \C{U} \times \C{P} \to \C{U} $
is some parameter dependent mapping like \refeq{eq:I-p} 
in \refS{intro}, which is well-posed in the sense that
for each $p\in\C{P}$ it has a unique solution $r(p)$ satisfying
\begin{equation}   \label{eq:II}
F(r(p),p) = 0,
\end{equation}
implicitly defining the function $r(p)$ alluded to above.
The mapping $F$ is representative for the conditions which
$r(p)$ has to satisfy.
\ignore{
one of the objects alluded to in the introduction, where $\C{P}$
is some set, and $\C{V}$ for the sake of
simplicity is assumed as a separable Hilbert space with inner product
$\langle\cdot|\cdot\rangle_{\C{U}}$ (the meaning of the index $\C{U}$ will soon
become clear).
}
What we desire is a
simple representation / approximation of that function, which avoids solving
\refeq{eq:II} every time one wants to know $r(p)$ for a new $p\in\C{P}$,
i.e.\ a \emph{proxy-}  or \emph{surrogate model}.

Of course the relation \refeq{eq:II} or its possible source \refeq{eq:I-p}
not only defines $r(p)$, but they can be an important relation each
candidate has to satisfy as well as possible, and possibly other such
relations.  This is important, as a proxy-model will often be used
also in the sense of a model order reduction, so that the computed $r_c(p)$
will be an approximation.  Then the degree to which a relation like \refeq{eq:II}
is satisfied can be the basis for estimating how good a particular approximation
$r_c(p)$ is.

One relatively well-known way when dealing with random models
\citep{segal58-TAMS, LGross1962, gelfand64-vol4, segalNonlin1969, kreeSoize86}
turns the problem into one of consideration and ultimately of
approximation of a linear mapping.   The details in the simplest case
are as follows.

\subsection{Associated linear map}   \label{SS:ass-lin-map}
Assume without significant loss of generality that
\ignore{
\begin{equation}    \label{eq:III}
\C{U} = \spn r(\C{P}) = \spn \im r \subseteq \C{V}
\end{equation}
}
$\spn r(\C{P}) = \spn \im r \subseteq \C{U} $,
the subspace of $\C{U}$ which is spanned by all the
vectors $\{r(p)|\; p\in\C{P}\}$,  is dense in $\C{U}$.

\begin{defi} \label{D:ass-lin}
Then to each such function $r:\C{P}\to \C{U}$ one may associate a linear map
\begin{equation}    \label{eq:IV}
R: \C{U} \ni u \mapsto \bkt{r(\cdot)}{u}_{\C{U}} \in \D{R}^{\C{P}},
\end{equation}
where $\D{R}^{\C{P}}$ is the space of real valued functions on $\C{P}$ and
$\bkt{r(\cdot)}{u}_{\C{U}}$ is the real valued map on $\C{P}$ given by
$\C{P}\ni p \mapsto \bkt{r(p)}{u}_{\C{U}}\in\D{R}$.
\end{defi}

\begin{lem} \label{L:inj}
By construction, $R$ restricted to $\spn\im r = \spn r(\C{P})$ is injective,
and hence has an inverse on its restricted range $\tilde{\C{R}}:=R(\spn \im r)$.
\end{lem}
\begin{proof}
Assume that for $u\in \im r=r(\C{P})$, $u\ne 0$, it holds that $R u = 0$.
This means that $\exists p_1\in\C{P}$ such that $u=r(p_1)$, and
$(Ru)(p)=\bkt{r(p)}{u}_{\C{U}}= \bkt{r(p)}{r(p_1)}_{\C{U}}=0$ for all $p\in\C{P}$.
Taking $p=p_1$, this means that $\bkt{r(p_1)}{u}_{\C{U}}=\bkt{r(p_1)}{r(p_1)}_{\C{U}}=
\nd{r(p_1)}_{\C{U}}^2=\nd{u}_{\C{U}}^2=0$.
This can only hold for $u = 0$, contradicting the assumption $u\ne 0$,
and so $R$ is injective on $\im r$ and by linearity
also on $\spn \im r$.  It is obviously also surjective from  $\spn \im r$ to $\tilde{\C{R}}$,
therefore \emph{bijective},  hence has an inverse $R^{-1}$ on $\tilde{\C{R}}$.
\end{proof}

\begin{defi}  \label{D:R-ip}
This may be used to define an inner product on $\tilde{\C{R}}$ as
\begin{equation}     \label{eq:V}
\forall \phi, \psi \in \tilde{\C{R}} \quad \langle \phi | \psi \rangle_{\C{R}} :=
     \langle R^{-1} \phi | R^{-1} \psi \rangle_{\C{U}},
\end{equation}
and to denote the completion of $\tilde{\C{R}}$ with this inner product by $\C{R}$.
\end{defi}

From \refL{inj} and \refD{R-ip} one immediately obtains

\begin{prop}  \label{P:R-unitary}
It is obvious from \refeq{eq:V} that $R^{-1}$ is a bijective isometry between $\spn \im r$
and $\tilde{\C{R}}$, hence continuous, and the same holds also $R$.  Hence extended by 
continuity to the completion Hilbert spaces, $R$ and $R^*=R^{-1}$ are
\emph{unitary} maps.
\end{prop}

Up to now, no structure on the set $\C{P}$ has been assumed, whereas on $\C{U}$
the inner product is assumed to measure what is important for the state 
$r(p) \in \C{U}$, i.e.\ vectors with large norm are considered important.
This is carried via the map $R$ defined in \refeq{eq:IV} onto
the space of scalar functions $\C{R}$ on the set $\C{P}$, and the inner product
there measures essentially the same thing as the one on $\C{U}$.
The only thing that changes is that now one does not have to work with
the space $\C{U}$, as everything is mirrored by the real functions
$\phi \in \C{R}$, which may be seen as a `problem-oriented co-ordinate system' on $\C{P}$.

\subsection{Reproducing kernel Hilbert space}   \label{SS:RKHS}
Given the maps $r:\C{P}\to\C{U}$ and $R:\C{U}\to\C{R}$, one may define the
\emph{reproducing kernel} \citep{berlinet, Janson1997}:

\begin{defi}  \label{D:RK} 
The reproducing kernel associated with $r:\C{P}\to\C{U}$ and $R:\C{U}\to\C{R}$
is $\vkappa \in \D{R}^{\C{P} \times\C{P}}$ is given by:
\begin{equation}     \label{eq:VI}
\C{P} \times\C{P} \ni (p_1, p_2) \mapsto
      \vkappa(p_1, p_2) := \bkt{r(p_1)}{r(p_2)}_{\C{U}} \in \D{R}.
\end{equation}
\end{defi}

It is straightforward to verify that:

\begin{thm}  \label{T:rep-prop}
For all $p\in\C{P}$: $\vkappa(p,\cdot)\in\tilde{\C{R}}\subseteq\C{R}$,
and $\spn \{ \vkappa(p,\cdot)\;\mid\; p\in\C{P} \}=\tilde{\C{R}}$, i.e.\
the kernel $\vkappa$ generates the space $\C{R}$.

The point evaluation functional $\updelta_p$ is a continuous map
on $\C{R}$, given by the inner product with the reproducing kernel:
\begin{equation}     \label{eq:VII}
 \updelta_p : \C{R} \ni \phi \mapsto \updelta_p(\phi) =
 \ip{\updelta_p}{\phi}_{\C{R}^* \times \C{R}}
     := \phi(p) = \bkt{\vkappa(p,\cdot)}{\phi}_{\C{R}} \in \D{R}.
\end{equation}
This reproduction of $\phi$ leads to the name \emph{reproducing kernel}.
\end{thm}
\begin{proof}
The first statement is due to the fact that 
$\vkappa(p,\cdot)=\ip{r(p)}{r(\cdot)}_{\C{U}} =Rr(p)(\cdot)$.
For the reproducing property,
similarly as in \refL{inj}, we take $\phi\in \tilde{\C{R}}$, i.e.\ $\exists u\in\C{U}$ with
$\phi(\cdot)=\bkt{r(\cdot)}{u}_{\C{U}} = R u(\cdot)$, and then extend by continuity to $\C{R}$.  
It holds for all $p\in\C{P}$:
\[ \updelta_p(\phi) = \bkt{\vkappa(p,\cdot)}{\phi}_{\C{R}} =
\bkt{R r(p)}{R u}_{\C{R}} = \bkt{r(p)}{u}_{\C{U}} = Ru(p) =\phi(p) ,\]
which is the reproducing property.  As $\updelta_p$ is defined via the inner
product, it is obviously continuous in $\phi$, hence this extends to the
closure of $\tilde{\C{R}}$, which is $\C{R}$.
\end{proof}

Hilbert spaces with such a reproducing kernel are called a \emph{reproducing kernel
Hilbert space} (RKHS)  \citep{berlinet, Janson1997}.
In other settings like classification or machine learning with support
vector machines, where $p \in \C{P}$ has to be classified as belonging
to a certain subsets of $\C{P}$, one can use such a map $r:\C{P}\to\C{U}$,
the so-called \emph{feature map}, implicitly through using an appropriate kernel.
This is then referred to as the `kernel trick', and classification may be
achieved by mapping these subsets with $r$ into $\C{U}$ and separating them
with hyperplanes---a linear classifier.  Observe also that the set
$\C{P}$ is embedded in $\C{R}$ via the correspondence
$\C{P}\ni p \mapsto \vkappa(p,\cdot) \in \C{R}$, which is the Riesz-representation
in $\C{R}$ of the continuous linear functional $\updelta_p \in \C{R}^*$.

Now we have a Hilbert space $\C{R}$ of real-valued
functions on $\C{P}$ and a linear surjective map
$R^{-1}:\C{R}\to\C{U}$ which can be used for representation.
In fact, as $\C{U}$ was assumed separable, so is the isomorphic space $\C{R}$,
one may now choose a basis $\{\vphi_m\}_{m\in\D{N}}$
in $\C{R}$, which may be assumed to be a complete orthonormal system (CONS).

\begin{coro} \label{C:RKHS-decomp}
With the CONS $\{y_m \;|\; y_m= R^{-1} \vphi_m = R^*\vphi_m\}_{m\in\D{N}}$ in $\C{U}$,
the operator $R$, and its adjoint or inverse $R^*=R^{-1}$, and the parametric
element $r(p)$ become
\begin{equation} \label{eq:VII0}
  R = \sum_m \vphi_m \otimes y_m;  \quad R^*=R^{-1} = \sum_m y_m \otimes \vphi_m;
  \quad r(p) = \sum_m \vphi_m(p) y_m = \sum_m R^* \vphi_m .
\end{equation}
These decompositions may also be seen as the \emph{singular value decompositions}
of the maps $R$ and $R^*=R^{-1}$, and are akin to the \emph{\KL{} decomposition} of $r(p)$.
\end{coro}
\begin{proof}
As $R$ is unitary, its singular values are all equal to unity, and any CONS such as
$\{\vphi_m\}_m$ is a set of \emph{right singular vectors}, giving the SVD of $R$
and hence of $R^*=R^{-1}$.    For any
$m, n\in\D{N}$ one has from the CONS property of $\{\vphi_m\}_{m}$ that
\[ \bkt{y_m}{y_n}_{\C{U}}=\bkt{R y_m}{R y_n}_{\C{R}}=\bkt{\vphi_m}{\vphi_n}_{\C{R}}
=\updelta_{m,n} , \]
and hence for any $p\in\C{P}$ and any $n\in\D{N}$:
$\vphi_n(p) =  \bkt{r(p)}{y_n}_{\C{U}} = (R y_n)(p)$,
due to \refeq{eq:IV} and the definition of the CONS $\{y_m\}$.  The last in
\refeq{eq:VII0} follows from definition of $R^*$.
\end{proof}

Observe that the relations \refeq{eq:VII0} exhibit the tensorial nature of the
representation mapping.  Looking at the \KL{} representation of $r(p)$,
one may see two things.  One is that this tensorial decomposition is a \emph{separated}
representation, as the $p$-dependence and the vector space have been separated.
The other is the observation
that \emph{model reductions} may be achieved by choosing only subspaces of $\C{R}$,
i.e.\ a---typically finite---subset of $\{\vphi_m\}_{m}$.  A good reduced order
model (ROM) is hence a representation where one only needs a few terms for a good
approximation.  Which subsets give a good
approximation will be addressed in the next \refS{correlat}.
Furthermore, the representation of $r(p)$ in \refeq{eq:VII0} is \emph{linear}
in the new `parameters' $\vphi_m$.  This means that by choosing the \emph{`co-ordinates'}
$\vphi_m$ on $\C{P}$, i.e.\ transforming $\C{P}\ni p \mapsto
(\vphi_1(p),\dots,\vphi_m(p),\dots)\in\D{R}^{\D{N}}$,
one obtains a \emph{linear / affine} representation on $\D{R}^{\D{N}}$.

%
%
%
%
%
%
%
%
%
%


%

\section{Correlation}  \label{S:correlat}
As already alluded to at the end of \refSS{ass-lin-map}, the RKHS construction $\C{R}$
with the inner product $\bkt{\cdot}{\cdot}_{\C{R}}$ just mirrors or reproduces the
inner product structure on the original space $(\C{U},\bkt{\cdot}{\cdot}_{\C{U}})$
on the RKHS space of real-valued functions $\C{R}$.
Up to now there is no way of telling what is important in the parameter set $\C{P}$.
Closely connected to this question is the one which subset of functions to
choose for model reduction.  Unfortunately,
up to now we have no indication which subset may be particularly good.
For this one needs additional information, a topic which will be taken up now.

As a way of indicating what is important on the set $\C{P}$, assume that
there is another inner product $\bkt{\cdot}{\cdot}_{\C{Q}}$ for
scalar functions $\phi \in \D{R}^{\C{P}}$, and denote the Hilbert
space of functions with that inner product by $\C{Q}$.
Abusing the notation a bit, we denote the map $R:\C{U}\to\C{Q}$, defined as in
\refeq{eq:IV} but with range $\C{Q}$, still by $R$.  Generally one would also assume
that the subspace $\dom R = \{ u\in\C{U} \mid \nd{Ru}_{\C{Q}}< \infty\}$ is, if
not the whole space $\C{U}$, at least dense in $\C{U}$.  Additionally, one would
assume that the densely defined operator $R$ is closed.
For simplicity assume here that $R$ is
defined on the whole space and hence continuous.  Furthermore, assume that the map 
$R:\C{U}\to\C{Q}$ is still injective, i.e.\  for $\phi\in\C{R}$ one has
$\nd{\phi}_{\C{R}}\ne 0 \Rightarrow \nd{\phi}_{\C{Q}}\ne 0$, and that $R$ is closed.
Without loss of generality we assume then that $R$ is surjective---by restricting ourselves
to the closed Hilbert subspace $R(\C{U})$ which we may call again $\C{Q}$.

\begin{defi}[Correlation]  \label{D:corr}
With this, one may define
\citep{kreeSoize86} a densely defined map $C$ in $\C{U}$ through the bilinear form 
\begin{equation}   \label{eq:IX}
   \forall u, v \in \C{U}:\quad \bkt{Cu}{v}_{\C{U}} := \bkt{Ru}{Rv}_{\C{Q}} .
\end{equation}
The map $C$, which may also be written as $C=R^* R$, may be called the 
\emph{`correlation'} operator.  By construction it is self-adjoint
and positive.  In case $R$ is defined on the whole space and hence continuous,
so is $C$.
\end{defi}

The last statements are standard results from the theory of linear operators
\citep{yosida-fa-1980}.  Observe that in contrast to \refSS{RKHS} the adjoint
is now different from the inverse as normally $R$ is not unitary,
i.e.\ the adjoint is of the map $R:\C{U}\to\C{Q}$
w.r.t.\ the $\C{Q}$-inner product $\bkt{\cdot}{\cdot}_{\C{Q}}$.
Often the inner product $\bkt{\cdot}{\cdot}_{\C{Q}}$ comes from a
measure $\vpi$ on $\C{P}$, so that for two measurable scalar functions
$\phi$ and $\psi$ on $\C{P}$ one has
\[
  \bkt{\phi}{\psi}_{\C{Q}} := \int_{\C{P}} \phi(p) \psi(p) \; \vpi(\di p),
\]
where the space $\C{Q}$ may then be taken as $\C{Q}:=\Lp_2(\C{P},\vpi)$;
or more generally with some kernel $\beta(p_1,p_2)$
\[
  \bkt{\phi}{\psi}_{\C{Q}} :=
  \int_{\C{P}\times\C{P}} \phi(p_1) \beta(p_1,p_2) \psi(p_2) \; \vpi(\di p_1) \vpi(\di p_2)
  = \ip{\beta}{\phi\otimes\psi}.
\]
One important sub-class of
such situations is when $\vpi$ is a probability measure on $\C{P}$, i.e.\
$\vpi(\C{P}) = 1$.  This is where the name `correlation' is borrowed from.
In the first case
\[
C = R^* R = \int_{\C{P}} r(p) \otimes r(p) \; \vpi(\di p).
\]
Often the set $\C{P}$ has more structure, like being in a topological space,
a differentiable (Riemann) manifold, or a Lie group, which then may induce
the choice of $\sigma$-algebra or measure.

\subsection{Spectral decomposition}   \label{SS:spec-dec}
Before, in \refSS{RKHS} it was the factorisation of $C= R^*R$ which allowed
the RKHS representation in \refeq{eq:VII0}.   For other representations, one needs other
factorisations.  Most common is to use the spectral
decomposition (e.g.~\citep{gelfand64-vol3, gelfand64-vol4, Segal1978,
reedSimon-vol1, reedSimon-vol2, DautrayLions3}) of $C$ to achieve such a factorisation.
In case the correlation were defined on a finite-dimensional space,
represented as a matrix $\vek{C}$, the eigenvalue problem---where $\lambda\in\D{C}$ is an
eigenvalue iff $\vek{C}-\lambda\vek{I}$ is not invertible---and eigen-decomposition
would be written with eigenvectors $\vek{v}_m$ and eigenvalues $\lambda_m$ as
\begin{equation}  \label{eq:ev-fd}
  \vek{C}\vek{v}_m = \lambda_m \vek{v}_m,\quad \vek{C}
  =\sum_m \lambda_m \vek{v}_m \vek{v}_m^{\ops{T}}=\sum_m \lambda_m \vek{v}_m \otimes\vek{v}_m
  = \sum_m \lambda_m \upDelta\vek{E}_m = \vek{V}\vek{\Lambda}\vek{V}^{\ops{T}} .
\end{equation}
As $\vek{C}$ is self-adjoint and positive, this implies $\lambda\in\D{R}$ and $\lambda_m\ge 0$.
The set of all eigenvalues $\sigma(\vek{C}):=\{\lambda_m\}_m\subset\D{C}$
is called the \emph{spectrum} of $\vek{C}$.  Here we assume the ordering
$0\le\lambda_1\dots\le\lambda_m$, each eigenvalue counted with appropriate multiplicity.
The $\vek{v}_m$ are normalised eigenvectors, and are mutually orthogonal 
($\vek{v}_m^{\ops{T}}\vek{v}_n=\updelta_{m,n}$).
The first two decompositions---which are only different notations---are into weighted sums
of simple tensor products of orthonormal vectors, or one-dimensional orthogonal projections
$\upDelta \vek{E}_m := \vek{v}_m \vek{v}_m^{\ops{T}} = \vek{v}_m \otimes\vek{v}_m$,
which define the spectral resolution $\vek{E}_m := \sum_{k\le m} \upDelta \vek{E}_k$.
The $\vek{E}_m$ are hence the orthogonal projections onto the subspaces
$\spn\{\vek{v}_k \mid k \le m \}$.

The columns of $\vek{V}=[\vek{v}_1,\dots,\vek{v}_m,\dots]$ are the normalised eigenvectors,
so that $\vek{V}$ is unitary resp.\ orthogonal,
and $\vek{\Lambda} = \diag(\lambda_m)$ is a diagonal matrix \citep{strang},
a `multiplication' operator, as for $\vek{\Lambda}\vek{u}=\vek{w}$, each component
$u_m$ of $\vek{u}=[u_1,\dots, u_m, \dots]^{\ops{T}}$ is just multiplied by $\lambda_m$:
$w_m = \lambda_m u_m$.  The last decomposition in \refeq{eq:ev-fd} hence means that
$\vek{C}$ is unitarily equivalent to a multiplication operator by a diagonal matrix
with real non-negative entries.

In contrast, on infinite dimensional Hilbert spaces the decompositions in \refeq{eq:ev-fd}
are materially different formulations of the spectral theorem for self-adjoint
operators, e.g.~\citep{gelfand64-vol3, Segal1978, reedSimon-vol1, reedSimon-vol2, DautrayLions3}.
A number $\lambda\in\D{C}$ is in the
spectrum $\sigma(C)$ iff $C-\lambda I$ is not invertible as a \emph{continuous operator}.
But now there may be spectral values $\lambda\in\sigma(C)$ which are \emph{not}
eigenvalues---this has to do with the possibility of a \emph{continuous spectrum}---and
the sums in \refeq{eq:ev-fd} have to become integrals.  Probably best known is the generalisation
of the second last form in \refeq{eq:ev-fd} ($\vek{C}=\sum_m \lambda_m \upDelta\vek{E}_m$),
namely \citep{gelfand64-vol3, Segal1978, reedSimon-vol1, reedSimon-vol1, DautrayLions3}:

\begin{thm}[First spectral theorem]     \label{T:1st-spec}
The self-adjoint and positive operator $C:\C{U}\to\C{U}$, where $C=R^* R$,
may be decomposed into an integral of orthogonal projections $E_\lambda$,
\begin{equation}   \label{eq:XII}
C = \int_0^\infty \lambda \; \di E_\lambda = \int_{\sigma(C)} \lambda \; \di E_\lambda.
\end{equation}
Here $E_\lambda$ is the corresponding projection-valued spectral measure corresponding to
$\vek{E}_m$ in \refeq{eq:ev-fd},
with a non-negative spectrum $\sigma(C) \subseteq \D{R}_+$.
\end{thm}
Observe that the factorised form $C=R^* R$ is actually equivalent to the statement that
$C$ is self-adjoint and positive.

For the sake of brevity and simplicity
of exposition let us assume that $C$ has a pure point spectrum $\sigma_p(C) = \sigma(C)$,
i.e.\ all $\lambda_m \in \sigma_p(C)$ are eigenvalues with eigenvector
$v_m$, $C v_m = \lambda_m v_m, m\in\D{N}$,
each eigenvalue repeated with appropriate finite multiplicity.
In this case \refeq{eq:XII}
becomes just a sum, and may be written with the CONS of unit-$\C{U}$-norm
eigenvectors $\{v_m\}_m \subset \C{U}$.  Here we assume the opposite
ordering of the $\lambda_m$ as before in \refeq{eq:ev-fd}, namely
$\lambda_1\ge\dots\ge\lambda_m\dots\ge 0$, and set
\begin{equation}  \label{eq:def-E_m}
   E_0 := I,\quad E_m := \sum_{k>m} v_k\otimes v_k; \quad \text{and for } m\ge 1: \;
      \upDelta E_{\lambda_m} := E_{m-1} - E_{m}.
\end{equation}
The spectral projection-valued measure $\di E_\lambda$ in \refeq{eq:XII}
becomes a point measure $\di E_\lambda = \sum_m \updelta_{\lambda_m}
\upDelta E_{\lambda_m}$, where $\updelta_{\lambda_m}$ is the \emph{Dirac}-$\updelta$.
For the second part of the following \refT{1st-spec-rep}, also assume that
the correlation $C$ is a \emph{trace class} or \emph{nuclear} operator, which means
that the trace is finite ($\tr C = \sum_m \lambda_m < \infty$), and
$C$ is then necessarily also a Hilbert-Schmidt and a compact operator.

\begin{thm}[First spectral representation and \KL{} expansion]  \label{T:1st-spec-rep}
The spectral decomposition of \refT{1st-spec}, \refeq{eq:XII} becomes
\begin{equation}   \label{eq:XIII}
C  = \sum_m \lambda_m (v_m \otimes v_m) = \sum_m \lambda_m \upDelta E_{\lambda_m}.
\end{equation}
Define a new CONS $\{s_m\}_m$ in $\C{Q}$:
 $\lambda_m^{1/2} s_m := R v_m$, to obtain the corresponding
\emph{singular value decomposition} (SVD) of $R$ and $R^*$.
The set $\vsigma(R)=\{\lambda_m^{1/2}\} = \sqrt{\sigma(C)}\subset \D{R}_+$
are the \emph{singular values} of $R$ and $R^*$.
\begin{equation}   \label{eq:XIV}
R = \sum_m \lambda_m^\frac{1}{2} (s_m \otimes v_m); \quad
R^* = \sum_m \lambda_m^\frac{1}{2} (v_m \otimes s_m); \quad
r(p) =  \sum_m \lambda_m^\frac{1}{2} \, s_m(p) v_m = \sum_m R^* s_m.
\end{equation}
\ignore{
\begin{eqnarray}  
    \label{eq:XIV} R = \sum_m \lambda_m^{1/2} s_m \otimes v_m, \quad
    R u &=& \sum_m \lambda_m^{1/2}\langle v_m | u \rangle_{\C{U}} s_m,\\
    \label{eq:XIVa}  R^* = \sum_m \lambda_m^{1/2} v_m \otimes s_m, \quad
    R^* \phi &=& \sum_m \lambda_m^{1/2}\langle s_m | \phi \rangle_{\C{W}} v_m,\\
    \label{eq:XIVb}  R^{-1} = \sum_m \lambda_m^{-1/2} v_m \otimes s_m, \quad
    R^{-1} \phi &=& \sum_m \lambda_m^{-1/2}\langle s_m | \phi \rangle_{\C{R}}v_m.
\end{eqnarray}
From this follows
\begin{equation}   \label{eq:XV}
  r(p) =  \sum_m \lambda_m^{1/2} \, s_m(p) v_m .
\end{equation}
}
The last relation is
the so-called \emph{\KL{} expansion} or \emph{proper orthogonal decomposition} (POD).
If in that relation the sum is \emph{truncated} at $n\in\D{N}$, i.e.\
\begin{equation}   \label{eq:best-n-term}
r(p) \approx r_n(p) =  \sum_{m=1}^n \lambda_m^\frac{1}{2} \, s_m(p) v_m = \sum_{m=1}^n R^* s_m(p),
\end{equation}
we obtain the \emph{best $n$-term approximation} to $r(p)$ in the norm of $\C{U}$.
\end{thm}
\begin{proof}
The spectral decompositions \refeq{eq:XIII}---analogues of the first three in
\refeq{eq:ev-fd}---are a consequence of the fact that for a point spectrum
the projection-valued measure $\di E_\lambda$ in \refeq{eq:XII} becomes a
discrete projection-valued measure $\upDelta E_{\lambda_m}$.

That the system $\{s_m\}_m$ is a CONS follows from
\[ 
\bkt{C v_m}{v_n}_{\C{U}} = \lambda_m \updelta_{m,n} = \bkt{R v_m}{R v_n}_{\C{Q}}
= \bkt{\lambda_m^{1/2} s_m}{\lambda_n^{1/2} s_n}_{\C{Q}} = \lambda_m \bkt{s_m}{s_n}_{\C{Q}} .
\]
The representations in \refeq{eq:XIV} are shown in the same way as in \refC{RKHS-decomp}.
It still remains to show that the function $p\mapsto r(p)$ defined in \refeq{eq:XV} is
in $\C{U}\otimes\C{Q}$.  For that, using the orthonormality of $\{s_m\}_m$ and $\{v_m\}_m$,
and the nuclearity of $C$, one computes
\[
\nd{r}_{\C{U}\otimes\C{Q}}^2 = \bkt{r}{r}_{\C{U}\otimes\C{Q}} = \sum_{m,n}
   \sqrt{\lambda_m \lambda_n}\,\bkt{s_m}{s_n}_{\C{Q}}\bkt{v_m}{v_n}_{\C{U}}  = \sum_{m,n}
   \sqrt{\lambda_m \lambda_n}\,\updelta_{m,n} \updelta_{m,n} = \sum_m \lambda_m < \infty.
\]
The statement about the best-$n$-term approximation follows
from the well-known optimality \citep{strang, Janson1997} of the SVD. 
\end{proof}

Observe that, similarly to \refeq{eq:VII0}, $r$ is linear in the $s_m$.
This means that by choosing the `co-ordinate transformation'
$\C{P}\ni p\mapsto (s_1(p),\dots,s_m(p),\dots)\in\D{R}^{\D{N}}$ one obtains
a \emph{linear / affine} representation where the first co-ordinates are
the most important ones, i.e.\ they catch most of the variability in that the
best-$n$-term approximation in the norm $\nd{\cdot}_C$ requires only the
first $n$ co-ordinate functions $\{s_m\}_{m\le n}$.  This is one possible
criterion on how to build good reduced order models $r_n(p)$, i.e.\ how to
choose a good subspace for approximation.

Note that in case $\C{P}$ is a probability space, the condition that $C$ be
a trace class or nuclear operator is also a necessary condition that $r$ have
finite variance and that the distribution of $r$ be a \emph{probability measure}
on $\C{U}$.  When stating other series representations in the sequel,
it will always be assumed that this condition of nuclearity is satisfied.
Hence the definition of models via linear maps is much more general
and allows one to consider \emph{generalised} resp.\ \emph{weak}, or in
some way ideal representations
\citep{segal58-TAMS, LGross1962, gelfand64-vol4, segalNonlin1969, kreeSoize86}.

\subsection{Singular value decomposition}   \label{SS:SVD}
To treat the analogues of the first two decompositions of $C$ in \refeq{eq:ev-fd}
in the case where $C$ has also a continuous spectrum directly requires
technical tools such as Gel'fand triplets (rigged Hilbert spaces), direct integrals
of Hilbert spaces \citep{gelfand64-vol3, gelfand64-vol4, DautrayLions3},
and generalised eigenvectors, which are beyond the scope of this short note.  
This also applies to representations
which go beyond the case of a nuclear correlation, and typically become some
kind of integral transform.  We contend ourselves with an alternative, and materially
stronger, formulation of the spectral decomposition than \refeq{eq:XII}
\citep{gelfand64-vol3, Segal1978, reedSimon-vol1, DautrayLions3},
the analogue of the last decomposition in \refeq{eq:ev-fd}
($\vek{C}=\vek{V}\vek{\Lambda}\vek{V}^{\ops{T}}$) which will lead us to the singular value
decomposition, an analogue of \refeq{eq:XIV}.  The results in this \refSS{SVD} do
not require $C$ to be nuclear, nor do they require $C$ or $R$ to be continuous.

\ignore{
It is of course also possible to directly use the unitary equivalence of
$L_2(\sigma(C))$ with $\C{U}$ and define for an appropriate $z\in\C{U}$
\begin{equation}  \label{eq:XVIaa}
\check{R}: L_2(\sigma(C)) \ni f \mapsto \int_0^\infty \lambda^{1/2}
  f(\lambda)\; \di E_\lambda z \in \C{U}.
\end{equation}
For a pure point spectrum \refeq{eq:XVIaa} reduces to
\begin{eqnarray} \label{eq:XVIab}
  \check{R}: f &\mapsto& \sum_m \lambda_m^{1/2} f(\lambda_m) v_m,\quad\text{or}\\
  \check{R} &=& \sum_m \lambda_m^{1/2} v_m \otimes \delta_{\lambda_m},
  \label{eq:XVIac}
\end{eqnarray}
where $\langle \delta_{\lambda_m}, f \rangle_{L_2(\sigma(C))} = f(\lambda_m)$,
again exhibiting the tensorial nature of the representation.
One finds that $\check{R} \check{R}^*=C$, and this is another
spectral representation.

An alternative, and materially stronger, formulation of the spectral decomposition
than \refeq{eq:XII} is  \citep{gelfand64-vol3, Segal1978, reedSimon-vol1, DautrayLions3} 
the analogue of the last decomposition in \refeq{eq:ev-fd}
($\vek{C}=\vek{V}\vek{\Lambda}\vek{V}^{\ops{T}}$):
}

\begin{thm}[Second spectral theorem]     \label{T:2nd-spec}
The densely defined, self-adjoint and positive operator
$C:\C{U}\to\C{U}$ is unitarily equivalent with a multiplication operator $M_{\mu}$,
\begin{equation}  \label{eq:ev-mult}
  C = V M_{\mu} V^*,
\end{equation}
where $V:\Lp_2(\C{T})\to\C{U}$ is unitary between some $\Lp_2(\C{T})$ on a measure
space $\C{T}$ and the Hilbert space $\C{U}$, and $M_{\mu}$ is a multiplication operator,
multiplying any $\psi\in\Lp_2(\C{T})$ with a real-valued
function $\mu$ --- $\Lp_2(\C{T})\ni\psi\mapsto \mu \psi\in\Lp_2(\C{T})$.
In case $C$ is bounded, ${\mu}\in\Lp_\infty(\C{T})$.
As $C$ is positive, $\mu(t)\ge 0$ for almost all $t\in\C{T}$,
and the essential range of $\mu$ is the spectrum of $C$.
As $M_{\mu}$ with a real valued non-negative $\mu$ is self-adjoint and positive,
one may define
\begin{equation}  \label{eq:ev-mult-sqrt}
  M_{\mu}^{1/2} := M_{\sqrt{\mu}}:\quad \Lp_2(\C{T})\ni\psi\mapsto \sqrt{\mu} \psi\in\Lp_2(\C{T}),
\end{equation}
from which one obtains the square-root of $C$ via its spectral decomposition
\begin{equation}  \label{eq:ev-mult-sqrt-2}
C^{1/2} = V  M_{\sqrt{\mu}} V^*.
\end{equation}
The factorisation corresponding to $C=R^* R$ in \refT{1st-spec} is here
(with $M_{\sqrt{\mu}}=M_{\sqrt{\mu}}^*$)
\begin{equation}  \label{eq:ev-mult-svd-2}
C= (V M_{\sqrt{\mu}})(V M_{\sqrt{\mu}})^* = (V M_{\sqrt{\mu}}I)(V M_{\sqrt{\mu}}I)^*. 
\end{equation}
\end{thm}
\begin{proof}
The statement about the unitary equivalence is a standard result 
\citep{gelfand64-vol3, Segal1978, reedSimon-vol1, DautrayLions3}
for self-adjoint operators, as well as the positivity of the multiplier $\mu$.
Computation of the square-root $M_{\sqrt{\mu}}$ is obvious, as $M_{\sqrt{\mu}}^2=M_{\mu}$;
on the other hand \refeq{eq:ev-mult-sqrt-2} is standard functional calculus of operators.
In the last \refeq{eq:ev-mult-svd-2}, it obviously holds that
$(V M_{\sqrt{\mu}})(V M_{\sqrt{\mu}})^*  = V M_{\sqrt{\mu}} M_{\sqrt{\mu}} V^* =  V M_{\mu} V^* $;
observe that $M_{\mu}^*=M_{\mu}$ and $M_{\sqrt{\mu}}^*= M_{\sqrt{\mu}}$, as $\mu$ is real.
\end{proof}

From this spectral decomposition follow decompositions of $R$ and some spectrally
connected factorisations of $C$:

\begin{coro}[Singular value decomposition and further factorisations] \label{C:SVD-R-fact}
The singular value decomposition (SVD) of $R$ is
\begin{equation}  \label{eq:ev-mult-svd}
R = U M_{\sqrt{\mu}} V^*,
\end{equation}
where $U:\Lp_2(\C{T})\to\C{Q}$ is a unitary operator,
$M_{\sqrt{\mu}}$ is from \refeq{eq:ev-mult-sqrt}, and the unitary $V$ from
\refeq{eq:ev-mult}.
Further decompositions of $C$ arising from \refT{2nd-spec} are $C=G^* G$ with
$G := I M_{\sqrt{\mu}} V^*$, and $C = (C^{1/2})^* C^{1/2} = C^{1/2} C^{1/2}$,
with the SVD of $C^{1/2}$ given by \refeq{eq:ev-mult-sqrt-2}.
\end{coro}
\begin{proof}
The SVD of $R$ in \refeq{eq:ev-mult-svd} is a standard result \citep{Segal1978, reedSimon-vol1},
and $U$ is unitary as $R$ was assumed surjective.  The decomposition with $G$
follows directly from \refeq{eq:ev-mult-svd-2}.   The last decomposition
$C = (C^{1/2})^* C^{1/2}$ follows from the fact that with $C$ also $C^{1/2}$
is self-adjoint, and as $C^{1/2}$ is also positive, its SVD is equal to its
spectral decomposition \refeq{eq:ev-mult-sqrt-2}.
\end{proof}

\subsection{Other factorisations and representations}   \label{SS:factor}
In the preceding \refSS{SVD} in \refC{SVD-R-fact} it was shown that there are
several ways to factorise $C=R^* R$.   Let us denote a general factorisation
by $C = B^* B$, where $B:\C{U}\to\C{H}$ is a map to a Hilbert space $\C{H}$
with all the properties demanded from $R$---see the beginning of this section.
Sometimes such a factor $B$ is called a \emph{square root} of $C$, but we
shall reserve that name for the \emph{unique} factorisation with the self-adjoint
factor $C^{1/2}$ from \refeq{eq:ev-mult-sqrt-2},
$C = (C^{1/2})^* C^{1/2}=C^{1/2} C^{1/2}$.
In some way, all such factorisations are equivalent:

\begin{thm}[Equivalence of factorisations]   \label{T:equi-fact}
Let  $C=B^*B$ with $B:\C{U}\to\C{H}$ be any factorisation satisfying the conditions
at the beginning of this section.
Any two such factorisations $B_1:\C{U}\to\C{H}_1$ and $B_2:\C{U}\to\C{H}_2$ 
with $C=B_1^*B_1=B_2^*B_2$ are
\emph{unitarily equivalent} in that there is a unitary map $X_{21}:\C{H}_1\to\C{H}_2$
such that $B_2 = X_{21} B_1$.  Equivalently, each such factorisation is unitarily
equivalent to $R$, i.e.\ for $C=B^* B$ there is a unitary $X:\C{H}\to\C{Q}$ such that
$R= X B$.
\end{thm}
\begin{proof}
Let $C=B_1^*B_1=B_2^*B_2$ be two such factorisations, each unitarily equivalent to
$R= X_1 B_1= X_2 B_2$.  As $X_2^*=X_2^{-1}$, it follows easily that
$B_2 = X_2^* X_1 B_1$, so $B_1$ and $B_2$ are unitarily equivalent with the
unitary $X_{21}:= X_2^* X_1$.

So what is left is to show that an arbitrary factorisation is equivalent to $R$.
From the SVD of $R$ in \refeq{eq:ev-mult-svd}, one sees easily that $R$ and $G$
in \refC{SVD-R-fact} are unitarily equivalent, as $R = U M_{\sqrt{\mu}} V^* 
= U (M_{\sqrt{\mu}} V^*) = U ( I M_{\sqrt{\mu}} V^*) = U G$.  Now let
$C=B^* B$ be an arbitrary factorisation with the required properties.  Then,
just as $R$ in \refC{SVD-R-fact}, the factor $B$ has a SVD \citep{Segal1978, reedSimon-vol1},
$B=W  M_{\sqrt{\mu}} V^*$, with $ M_{\sqrt{\mu}}$ and $V$ from \refC{SVD-R-fact},
and a unitary $W:\Lp_2(\C{T})\to\C{H}$.  Hence $B = W G$ or $G=W^* B$, and finally
$R = U G=U W^* B=X B$ with a unitary $X:=U W^*$.
\end{proof}

For finite dimensional spaces, a favourite choice for such a decomposition
of $C$ is the Cholesky factorisation $\vek{C} = \vek{L} \vek{L}^{\ops{T}}$,
where $B=\vek{L}^{\ops{T}}$.
Now let us go back to the situation described in \refT{1st-spec-rep}, where for the
sake of simplicity of exposition we assume that $C$ has a purely discrete spectrum
and a CONS of eigenvectors $\{v_m\}_m$ in $\C{U}$, and let us have
a look how the results up to now may be used to build new representations.
First transport the eigenvector CONS from $\C{U}$ to $\Lp_2(\C{T})$:

\begin{lem}   \label{L:eig-xi}
Setting for all $m\in\D{N}$: $\xi_m := V^* v_m$, the system $\{\xi_m\}_m$
is a CONS in $\Lp_2(\C{T})$, and $M_\mu \xi_m = \lambda_m \xi_m$,
i.e.\ the $\xi_m$ are an eigenvector CONS of $M_\mu = V^* C V$.
\end{lem}
\begin{proof}
Orthonormality and completeness are due to $V^*$ being unitary.  With
\refeq{eq:ev-mult} one computes
\[
  M_\mu \xi_m = V^* C V V^* v_m = \lambda_m V^* v_m = \lambda_m \xi_m,
\]
which shows the eigenvector property.
\end{proof}

\begin{prop}   \label{P:new-cons}
With the help of the CONS $\{s_m\}_m$ in $\C{Q}$ or $\{v_m\}_m$ in $\C{U}$,
define a CONS $\{h_m\}_m$ in $\C{H}$:  $\overline{\spn}\{h_m \mid m\in\D{N} \} = \C{H}$:
\begin{equation}  \label{eq:equi-cons}
   \forall m\in\D{N}: \quad h_m := B C^{-1} R^* s_m = B C^{-1/2} v_m.
\end{equation}
The CONS $\{h_m\}_m$ in $\C{H}$ is an eigenvector CONS of the operator
\begin{equation}  \label{eq:def-C-H}
   C_{\C{H}} := B B^*:\C{H}\to\C{H},
\end{equation}
\begin{equation}  \label{eq:eig-v-C-H}
   \forall m\in\D{N}: \quad C_{\C{H}} h_m := \lambda_m h_m.
\end{equation}
\end{prop}
\begin{proof}
The stated orthonormality of the $\{h_m\}_m$ is easily computed, as with
\refT{2nd-spec}, \refC{SVD-R-fact}, and the SVD of $B=W  M_{\sqrt{\mu}} V^*$
from \refT{equi-fact}, one obtains after a bit of computation $B C^{-1} R^* = W U^*$,
and $B C^{- 1/2} = W V^*$, hence $h_m = W U^* s_m$, and therefore orthonormality follows
from the unitarity of $W U^*$ and orthonormality of the $\{s_m\}_m$.
Completeness follows from the completeness of  $\{s_m\}_m$ and surjectivity of $B$.

Similarly to $h_m = W U^* s_m$ one obtains with  \refL{eig-xi}:
\[
  v_m =  V U^* s_m \Rightarrow h_m = W U^* s_m = W V^* (V U^* s_m) = W V^* v_m = W \xi_m .
\]
From this follows, again with \refL{eig-xi},
\[
   C_{\C{H}} h_m = (B B^*) W \xi_m = W M_\mu \xi_m = \lambda_m W \xi_m = \lambda_m h_m ,
\]
\end{proof}

One may see the statement in \refL{eig-xi} as a special case of Proposition~\ref{P:new-cons}
with $B = I M_{\sqrt{\mu}} V^*$, as then $\xi_m = V^* v_m = U^* s_m$.
Collecting, an immediate consequence is:

\begin{coro}  \label{C:other-eig}
One has the following equivalent eigensystems
\begin{itemize}
\item on $\C{U}$ with $C = R^* R$ --- $C v_m = \lambda_m v_m$,
       $v_m = V U^* s_m$ and $C = V M_\mu V^*$;
\item on $\C{H}$ with $C_{\C{H}} = B B^*$ --- $C_{\C{H}} h_m = \lambda_m h_m$,
       $h_m = W V^* v_m$ and $C_{\C{H}} = W V^* C V W^*$;
\item on $\C{Q}$ with $C_{\C{Q}} = R R^*$ --- $C_{\C{Q}} s_m = \lambda_m s_m$,
       $s_m = U V^* v_m$ and $C_{\C{Q}} = U V^* C V U^*$;
\item on $\Lp_2(\C{T})$ with $C_{\Lp_2(\C{T})} = M_\mu$ ---
       $C_{\Lp_2(\C{T})} \xi_m = \lambda_m \xi_m$, $\xi_m = V^* v_m$ 
       and $C_{\Lp_2(\C{T})} = V^* C V$.
\end{itemize}
\end{coro}
The last two statements are special cases of the second one with $\C{H}=\C{Q}$ and
$B = R$, resp.\ $\C{H}=\Lp_2(\C{T})$ with $B = M_{\sqrt{\mu}} V^*$.
Hence each factorisation $C=B^* B$ with $B:\C{U}\to\C{H}$ gives
a new equivalent eigensystem on $\C{H}$ for the operator $C_{\C{H}} = B B^*$.

From \refeq{eq:XIV} in \refT{1st-spec-rep} one has $r = \sum_m R^* s_m$.  This
in conjunction with another equivalent factorisation according to \refT{equi-fact}
immediately leads to new representations of $r(p)$, by replacing $R^* s_m$ in the
\KL{} expansion in \refT{1st-spec-rep} by the equivalent $B^* h_m$.

\begin{coro}[Representation from factorisation]  \label{C:fact-rep}
With a factorisation $C=B^* B$ and  CONS $\{h_m\}_m$ in $\C{H}$ as in
Proposition~\ref{P:new-cons}, one obtains
the following representation of $r(p)$:
\begin{equation}  \label{eq:fact-rep}
   r(p) = \sum_m B^* h_m= \sum_m V M_{\sqrt{\mu}}W^* h_m, \quad
   \text{in particular also} \quad 
     r(t) = \sum_m V M_{\sqrt{\mu}} \xi_m(t).
\end{equation}
\end{coro}

In the special case of a purely discrete spectrum we are dealing with here
it is possible to formulate analogues of the decompositions in \refC{RKHS-decomp}.
This is an analogue of \refT{1st-spec-rep} for
the general case $C=B^* B$  with $B:\C{U}\to\C{H}$: 

\begin{coro}  \label{C:fact-rep-tens}
With a factorisation $C=B^* B$ and CONS $\{v_m\}_m$ in $\C{U}$, CONS $\{h_m\}_m$
in $\C{H}$ as in Proposition~\ref{P:new-cons}, one obtains
the following tensor representations of the map $C_{\C{H}}= B B^*$:
\begin{equation}  \label{eq:C-tens-rep}
   C_{\C{H}} = \sum_m \lambda_m h_m\otimes h_m.
\end{equation}
Specifically, for $\C{H}=\C{Q}$ and $C_{\C{H}}=C_{\C{Q}}$,
\begin{equation}  \label{eq:C-tens-rep-Q}
   C_{\C{Q}} = \sum_m \lambda_m s_m\otimes s_m.
\end{equation}
The corresponding expansions of $B$ and its adjoint are:
\begin{equation}  \label{eq:B-tens-rep}
   B = \sum_m \lambda_m^{1/2}\, h_m\otimes v_m; \quad  
   \text{and} \quad B^* = \sum_m \lambda_m^{1/2}\, v_m\otimes h_m.
\end{equation}
In case the space $\C{H}$ is a function space like $\Lp_2(\C{T})$
on a set $\C{A}$, this results in the \KL{} expansions for the representation of $r(p)$:
\begin{equation}  \label{eq:r-fact-rep-tens}
   r(a) = \sum_m \lambda_m^{1/2}\, h_m(a) v_m, \quad
   \text{in particular also} \quad 
     r(t) = \sum_m \lambda_m^{1/2}\, \xi_m(t) v_m.
\end{equation}
\end{coro}

In this last \refeq{eq:r-fact-rep-tens} the function $r(p)$ has become a function
of the new parameter $a\in\C{A}$ or $t\in\C{T}$, having implicitly performed
a transformation $\C{P}\to\C{A}$ or $\C{P}\to\C{T}$.  The new parametrisation covers
the same range as $r(p)$ before.
As a summary of the analysis let us put everything together:

\begin{thm}[Equivalence of representation and factorisation] \label{T:}
A parametric mapping $r:\C{P}\to \C{U}$ into a Hilbert space $\C{U}$ with the conditions
stated at the beginning of \refSS{ass-lin-map} and this section induces a
linear map $R:\C{U}\to\C{Q}$, where $\C{Q}$ is a Hilbert space of functions
on $\C{P}$.  The \emph{reproducing kernel Hilbert space} is a special case of this.

Any other factorisation of the `correlation' $C=R^* R$ on $\C{U}$, like $C = B^* B$ with
a $B:\C{U}\to\C{H}$ into a Hilbert space $\C{H}$ with the same properties as $R$
is unitarily equivalent, i.e.\ there is a unitary $W:\C{Q}\to\C{H}$ such that
$B = W R$.  Any such factorisation induces a representation of $r$.  Especially
if $\C{H}$ is a space of functions on a set $\C{A}$, one obtains a representation
$r(a), (a\in\C{A})$, such that $r(\C{P}) = r(\C{A})$.

The associated `correlations' $C_{\C{Q}} = R R^*$ on $\C{Q}$ resp.\  $C_{\C{H}} = B B^*$
on $\C{H}$ have the same spectrum as $C$, and factorisations of $C_{\C{Q}}$ resp.\ $C_{\C{H}}$
induce new factorisations of $C$.
\end{thm}

%
%
%
%
%
%
%
%
%
%
%
%


%

\section{Kernel space}  \label{S:kernel}
In this section we take a closer look at the operator defined in \refeq{eq:def-C-H}
in Proposition~\ref{P:new-cons} especially for the case $\C{H}=\C{Q}$ and $B = R$, i.e.\
we analyse the operator $C_{\C{Q}} = R R^*$.  We shall restrict ourselves again
to the case of a pure point spectrum.  From \refC{other-eig} and
\refeq{eq:C-tens-rep-Q} in \refC{fact-rep-tens} one knows that in an abstract sense
\begin{equation}   \label{eq:C-Q-decs}
  C_{\C{Q}} = U V^* C V U^* = U M_\mu U^* = \sum_m \lambda_m s_m \otimes s_m .
\end{equation}
But the point is here to spell this out in more analytical detail especially
for the case when, as indicated already at the beginning of \refS{correlat},
the inner product on $\C{Q}$ is given by a measure $\vpi$ on $\C{P}$:
\begin{equation}   \label{eq:Q-ip}
\forall \vphi, \psi \in \C{Q}: \quad \bkt{\vphi}{\psi}_{\C{Q}} =
       \int_{\C{P}} \vphi(p) \psi(p)\, \vpi(\di p) .
\end{equation}

\subsection{Kernel spectral decomposition}   \label{SS:kernel-spec-dec}
Then $C$ is given  by 
\begin{equation}   \label{eq:C-tens-int}
  C = \int_{\C{P}} r(p)\otimes r(p)\, \vpi(\di p) ,
\end{equation}
and $C_{\C{Q}}$ is represented by the kernel
\begin{equation}   \label{eq:C-Q-kernel}
  \vkappa(p_1,p_2) = \bkt{r(p_1)}{r(p_2)}_{\C{U}}, 
\end{equation}
so that for all $\vphi, \psi \in \C{Q}$
\begin{equation}   \label{eq:XVIII}
   \bkt{C_{\C{Q}}\phi}{\psi}_{\C{Q}} = \bkt{R^* \vphi}{R^* \psi}_{\C{U}} = 
   \iint_{\C{P}\times\C{P}}  \vphi(p_1) \vkappa(p_1, p_2) \psi(p_2)\; \vpi(\di p_1) \vpi(\di p_2),
\end{equation}
i.e.\ $C_{\C{Q}}$ is a Fredholm integral operator
\begin{equation}   \label{eq:XVIII-int-op}
   (C_{\C{Q}} \psi)(p_1)  = \int_{\C{P}}  \vkappa(p_1, p_2) \psi(p_2)\; \vpi(\di p_2).
\end{equation}

The abstract eigenvalue problem described in \refC{other-eig} for the
operator $C_{\C{Q}}$, when taking into account the explicit description
\refeq{eq:XVIII-int-op}, is translated into finding an eigenfunction $s\in\C{Q}$
and eigenvalue $\lambda$ such that
\begin{equation}   \label{eq:XVIII-int-ev}
   (C_{\C{Q}} s)(p_1)  = \int_{\C{P}}  \vkappa(p_1, p_2) s(p_2)\; \vpi(\di p_2)
      = \lambda \, s(p_1) ,
\end{equation}
a \emph{Fredholm} integral equation \citep{courant_hilbert, atkinson97}.

\begin{prop}  \label{P:kernel-eig-fcts}
From \refC{other-eig} and \refeq{eq:C-Q-decs} one knows that the eigenfunctions
are $\{s_m\}_m \subset \C{Q}$, hence, in particular with the kernel $\vkappa$,
Mercer´s theorem \citep{courant_hilbert} gives
\begin{equation}   \label{eq:XVIIIb}
 \int_{\C{P}}  \vkappa(p_1, p_2) s_m(p_2)\; \vpi(\di p_2)
      = \lambda_m \, s_m(p_1); \quad
\vkappa(p_1, p_2) = \sum_m \lambda_m\, s_m(p_1) s_m(p_2),
\end{equation}
giving a decomposition of $\vkappa$, which is of course essentially
identical to \refeq{eq:C-Q-decs}.
\end{prop}

In \refS{correlat} the analysis was based to a large extent on factorisations of
the operator $C = R^* R$.   Similarly, now one looks at factorisations of
$C_{\C{Q}} = R R^*$.

One example situation which occurs quite frequently fits nicely here, rather than
in the later \refS{xmpls}, which is the case when $\C{P} = \D{R}^n$ with the usual
Lebesgue measure, and the kernel
is a convolution kernel, i.e.\  $\vkappa(p_1, p_2) = \kappa(p_1 - p_2)$.
This means that the kernel is invariant under arbitrary displacements or shifts
$z\in\D{R}^n$:  $\vkappa(p_1, p_2) = \vkappa(p_1 + z, p_2 + z) = \kappa(p_1 - p_2)$.
The eigenvalue equation \refeq{eq:XVIII-int-ev} becomes
\begin{equation}   \label{eq:XVIII-int-ev-z}
   (C_{\C{Q}} s)(p_1)  = \int_{\D{R}^n}  \kappa(p_1- p_2) s(p_2)\; \di p_2
      = \lambda \, s(p_1) .
\end{equation}
As is well known, the symmetry of $\vkappa$ implies now that the
function $\kappa$ has to be an \emph{even} function \citep{bracewell},
$\kappa(z) = \kappa(-z)$.

It is clear that this form of equation can be treated by Fourier analysis
\citep{courant_hilbert, atkinson97}; performing a Fourier transform on
\refeq{eq:XVIII-int-ev-z} and denoting transformed quantities by a hat, e.g. $\hat{s}$,
one obtains for all $\zeta \in \D{R}^n$
\begin{equation}   \label{eq:XVIII-int-ev-F}
   \widehat{(C_{\C{Q}} s)}(\zeta)  =  \hat{\kappa}(\zeta) \hat{s}(\zeta)
      = \lambda \, \hat{s}(\zeta) .
\end{equation}
In this representation, $\widehat{(C_{\C{Q}} s)}$ has become a multiplication
operator $M_{\hat{\kappa}}$ with the positive multiplier function
$\hat{\kappa}(\zeta)\ge 0$ on the domain $\zeta\in\D{R}^n$.
This is a concrete version of the case in the second spectral \refT{2nd-spec},
the multiplier function $\mu(t)$ in that theorem is $\hat{\kappa}(\zeta)$ here.
The unitary transformation which has effectively \emph{diagonalised} the integral
operator \refeq{eq:XVIII-int-ev} is the \emph{Fourier transform}, and the essential range
of $\hat{\kappa}$ is the spectrum.  This relates to the fact that the Fourier
transform of the correlation $\kappa$ --- or more precisely the covariance, but we do not
distinguish this here --- is usually called the spectrum, or more precisely the spectral
density.  In the terminology here the spectrum is the \emph{values} of $\hat{\kappa}$.
It is now also easy to see how a continuous spectrum appears: on an infinite domain
the integral operator \refeq{eq:XVIII-int-op} is typically not compact, and unless
$\kappa$ is almost-periodic the operator has no point spectrum.  The Fourier functions
are \emph{generalised} eigenfunctions \citep{gelfand64-vol3, gelfand64-vol4, DautrayLions3},
as they are \emph{not} in $\Lp_2(\D{R}^n)$.  We shall not dwell further on
this topic here.

If we denote the Fourier transform on $ \Lp_2(\D{R}^n)$ by 
\begin{equation}  \label{eq:fourier}
  F: f(p) \mapsto \hat{f}(\zeta)= (F f)(\zeta) = \int_{\D{R}^n} 
     \exp(-2 \uppi\, \ii\, p\cdot \zeta) f(p) \,\di p ,
\end{equation}
where $p\cdot \zeta$ is the Euclidean inner product in $\D{R}^n$,
then one may write this spectral decomposition and factorisation of $C_{\C{Q}}$
in this special case corresponding to \refC{fact-rep-tens} in the following way.

\begin{coro}  \label{C:rep-kern-spec}
The operator $C_{\C{Q}}=R R^*$ has
in the stationary case of \refeq{eq:XVIII-int-ev-z} the spectral decomposition
\begin{equation}  \label{eq:C-Q-spec-dec}
   C_{\C{Q}} = F^* M_{\hat{\kappa}} F.
\end{equation}
As $\hat{\kappa}(\zeta)\ge 0$, the square-root multiplier is given by
\begin{equation}  \label{eq:M-Q-spec-sqrt}
   M_{\hat{\kappa}}^{1/2} = M_{\sqrt{\hat{\kappa}}} .
\end{equation}
This induces the following factorisation of $C_{\C{Q}}=R R^*$:
\begin{equation}  \label{eq:C-Q-spec-fact}
   C_{\C{Q}} = (M_{\sqrt{\hat{\kappa}}}F)^* (M_{\sqrt{\hat{\kappa}}} F).
\end{equation} 
\end{coro}

From \refC{other-eig} one has $C_{\C{Q}} = U V^* C V U^*$, which gives further
\begin{equation}  \label{eq:fourier-dent-U-F}
  C_{\C{Q}} = U V^* C V U^* = U V^* V M_\mu V^* V U^* = U M_\mu U^*,
\end{equation}
meaning that essentially in this case $U = F^*$, the inverse Fourier transform.
This implies the well-known Fourier representations of stationary random functions.
Denoting the shift operator for $z\in\D{R}^n$ as $T_z: f(p)\mapsto f(p+z)$,
it is elementary that with $\eta_\zeta(p) := \exp(2 \uppi\, \ii\, p\cdot \zeta)$
\[
  T_z \eta_\zeta(p) = T_z \exp(2 \uppi\, \ii\, p\cdot \zeta) = \exp(2 \uppi\, \ii\, p\cdot z)\,
  \exp(2 \uppi\, \ii\, p\cdot \zeta) = \exp(2 \uppi\, \ii\, p\cdot z)\, \eta_\zeta(p),
\]
which says that the $\eta_\zeta(p)$ are `generalised' eigenfunctions of $T_z$
 \citep{gelfand64-vol4, DautrayLions3}.  They
are \emph{not} true eigenfunctions as they are not in $\Lp_2(\D{R}^n)$.

Shift-invariance means that $T_z\,C_{\C{Q}} = C_{\C{Q}}\, T_z$, i.e.\ the operators
commute.  This implies that  \citep{gelfand64-vol3, Segal1978, DautrayLions3} they have the same
spectral resolution, i.e.\ the same true and generalised eigenfunctions.  Both $T_z$
and $C_{\C{Q}}$ are effectively diagonalised by the Fourier transform $F$.
This particular case of covariance has been treated extensively in the literature
\citep{Karhunen1946, Karhunen1950,  Yaglom-62-04, gelfand64-vol4, Yaglom-68-I, Yaglom-68-II,
LoeveII, kreeSoize86, Matthies_encicl}.  As is well known, the functions
$\eta_\zeta(p)$ are formal or generalised
eigenfunctions of a shift-invariant operator as that one in \refeq{eq:XVIII-int-ev-z}
\citep{gelfand64-vol3, gelfand64-vol4, DautrayLions3}.
This results in the following \KL{} expansions, also known as spectral expansions,
for the formal representation of $r(p)$ in the case of a discrete spectrum
\refeq{eq:C-tens-rep-Q}:
\begin{equation}  \label{eq:r-fact-rep-C-Q-d} 
     r(p) = \sum_m \lambda_m^{1/2}\, \eta_{\lambda_m}(p) v_m,
\end{equation}
or, in the case of a continuous spectrum with generalised eigenvectors $v_\zeta$ of $C$,
formally,
\begin{equation}  \label{eq:r-fact-rep-C-Q-c} 
     r(p) = \int_{\D{R}^n} \sqrt{\hat{\kappa}(\zeta)}\, \eta_{\zeta}(p) v_\zeta \,\di \zeta =
 \int_{\D{R}^n}  \exp(2 \uppi\, \ii\, p\cdot \zeta)\,M_{\sqrt{\hat{\kappa}}}\, v_\zeta \,\di \zeta
     = F^* (M_{\sqrt{\hat{\kappa}}} v_\zeta),
\end{equation}
or, in the most general case, a combination of \refeq{eq:r-fact-rep-C-Q-d} and
\refeq{eq:r-fact-rep-C-Q-c}.  In that last \refeq{eq:r-fact-rep-C-Q-c}, 
the formal term $v_\zeta\,\di \zeta$
may be interpreted as a vector-valued measure $\tilde{v}(\di \zeta)$, in the case of
a \emph{random} process or field $r(p)$ on $p\in\D{R}^n$ it is called a \emph{stochastic}
measure.

\subsection{Kernel factorisations}   \label{SS:kernel-fact}
The concrete realisation of the operator $C_{\C{Q}}$ as an integral kernel
\refeq{eq:XVIII-int-op} opens the possibility to look for factorisations
in the concrete setting of integral transforms.

If, on the other hand, one has some other factorisation of the kernel,
for example on some measure space $(\C{X},\nu)$:
\begin{equation}  \label{eq:XIX}
  \vkappa(p_1,p_2) = \int_{\C{X}} g(p_1,x) g(p_2,x)\; \nu(\di x) = 
   \bkt{g(p_1,\cdot)}{g(p_2,\cdot)}_{\Lp_2(\C{X})},
\end{equation}
then the integral transform with kernel $g$ will play the r\^ole of a
factor as before the mappings $R$ or $B$.
Let us recall that in the context of RKHS, cf.\ \refSS{RKHS}, such a factorisation is
often used as a so-called \emph{feature map}.

\begin{defi} \label{D:int-trfm}
Define the integral transform $X: \Lp_2(\C{X}) \to \C{Q}$ with kernel $g$ as
\begin{equation}  \label{eq:XIX-int-trfm}
X: \xi \mapsto \int_{\C{X}} g(\cdot,x)\xi(x)\; \nu(\di x).
\end{equation}
\end{defi}

This results immediately in a new factorisation of $C_{\C{Q}}$ and a new representation:

\begin{coro}  \label{C:rep-kern-trfm}
The operator $C_{\C{Q}}=R R^*$ with decomposition \refeq{eq:C-Q-decs} has
the factorisation 
\begin{equation}  \label{eq:C-Q-fact-g-X}
   C_{\C{Q}} = X X^*.
\end{equation}
Defining the orthonormal system $\{\chi_m\}_m \subset \Lp_2(\C{X})$ by
\begin{equation}  \label{eq:XIX-k-trfm}
  \lambda_m^{1/2} \chi_m = X^* s_m; \quad \lambda_m^{1/2} \chi_m(x) =
  \int_{\C{P}} g(p,x) s_m(p)\; \vpi(\di p),
\end{equation}
this induces the following \KL{} representation of $r(x)$ on $\C{X}$:
\begin{equation}  \label{eq:r-C-Q-X-rep}
  r(x) = \sum_m \lambda_m^{1/2} \chi_m(x) v_m.
\end{equation} 
\end{coro}
\begin{proof}
  To prove \refeq{eq:C-Q-fact-g-X}, compute for any $\phi\in\C{Q}$ its adjoint transform
  $(X^* \phi)(x) = \int_{\C{P}} g(p,x) \phi(p)\; \vpi(\di p)$.  Now for all
  $\vphi, \psi \in\C{Q}$ it holds that
\begin{multline*}
   \bkt{X X^* \vphi}{\psi}_{\C{Q}} =  
   \bkt{X^* \vphi}{X^* \psi}_{\Lp_2(\C{X})} = \int_{\C{X}}  (X^* \vphi)(x)(X^* \psi)(x)
   \; \nu(\di x) = \\ \int_{\C{X}} \left(\int_{\C{P}} g(p_1,x)
   \vphi(p_1)\; \vpi(\di p_1) \right) \left( \int_{\C{P}} g(p_2,x) \psi(p_2)\; \vpi(\di p_2)
   \right) \, \nu(\di x) = \\
   \iint_{\C{P}\times\C{P}} \left( \int_{\C{X}} g(p_1,x) g(p_2,x)\; \nu(\di x) \right)
   \vphi(p_1) \psi(p_2)\; \vpi(\di p_1) \vpi(\di p_2) = \\
   \iint_{\C{P}\times\C{P}} \vkappa(p_1,p_2) \vphi(p_1) \psi(p_2)\; \vpi(\di p_1) \vpi(\di p_2)
   = \bkt{C_{\C{Q}} \vphi}{\psi}_{\C{Q}},
\end{multline*}
which is the bilinear form for \refeq{eq:C-Q-fact-g-X}.
The rest follows from \refC{fact-rep-tens} with $\C{H}=\Lp_2(\C{T})$ and
$B = X^* U V^*: \C{U}\to\Lp_2(\C{X})=\C{H}$, as from \refeq{eq:C-Q-decs}
\[
  X X^* = C_{\C{Q}} = U V^* C V U^*\quad \Rightarrow \quad  C = (V U^* X) (X^* U V^*) = B^* B .
\]
\end{proof}

The decomposition in Proposition~\ref{P:kernel-eig-fcts} may now also be seen in this light
by setting $\C{X}=\D{N}$ with \emph{counting} measure $\nu$, such that
$\Lp_2(\C{X}) = \ell_2$, and $X$-transformation kernel $g(p,m) := \lambda_m^{1/2} s_m(p)$.
Then \refeq{eq:XIX} becomes \refeq{eq:XVIIIb}, the concrete version of \refeq{eq:C-Q-decs}.

The result in
\refeq{eq:C-Q-spec-fact}, $C_{\C{Q}} = (M_{\sqrt{\hat{\kappa}}} F)^* (M_{\sqrt{\hat{\kappa}}} F)$
shows that the Fourier diagonalisation  in \refC{rep-kern-spec}
is a special case of such a kernel factorisation with $X := (M_{\sqrt{\hat{\kappa}}}\,F)^*$.
As $\kappa$ is the inverse Fourier transform of $\hat{\kappa}$,
\[
\kappa(p) = 
\int_{\D{R}^n}  \exp(2 \uppi\, \ii\, p\cdot \zeta)\,\hat{\kappa}(\zeta) \,\di \zeta ,
\]
remembering that with the Fourier transform one has to consider the \emph{complex}
space $\C{C}=\Lp_2(\D{R}^n,\D{C})$ with inner product
\[
\forall \vphi, \psi \in \C{C}: \quad \bkt{\vphi}{\psi}_{\C{C}} :=
   \int_{\D{R}^n} \bar{\vphi}(\zeta) \psi(\zeta) \,\di \zeta 
\] 
($\bar{\vphi}(\zeta)$ is the conjugate complex of $\vphi(\zeta)$), and by defining
the $X$-transformation kernel
$\gamma(p,\zeta) := \exp(2 \uppi\, \ii\, p\cdot \zeta)\,\sqrt{\hat{\kappa}(\zeta)}$,
one obtains the kernel factorisation
\begin{multline*}
   \vkappa(p_1,p_2) =  \kappa(p_1 - p_2) =
   \int_{\D{R}^n}  \exp(2 \uppi\, \ii\, (p_1 - p_2) \cdot \zeta)\, \hat{\kappa}(\zeta) 
   \,\di \zeta = \\ \int_{\D{R}^n}\left(\exp(-2\uppi\,\ii\,p_2\cdot\zeta)\,
   \sqrt{\hat{\kappa}(\zeta)}\right) \left(\exp(2\uppi\,\ii\,p_1\cdot\zeta)\,\sqrt{\hat{\kappa}
   (\zeta)}\right)\,\di \zeta = \bkt{\gamma(p_2,\cdot)}{\gamma(p_1,\cdot)}_{\C{C}}.
\end{multline*}

%
%
%
%
%
%
%
%
%
%
%
%


%

\section{Interpretations, decompositions, and reductions}  \label{S:xmpls}
After all the abstract deliberations it is now time to see some concrete examples,
which will show that the above description
is in many cases an abstract statement of already very familiar constructions.

An important example of these decompositions is when $\C{U}$ is also
a space of functions.  Imagine for example a scalar \emph{random field} $u(x,\omega)$,
where $x\in\C{X}\subset\D{R}^n$ is a spatial variable, and $\omega\in\Omega$ is
an elementary event in a probability space $\Omega$ with probability measure $\D{P}$.
Na\"ively, at each $x\in\C{X}$ there is a random variable (RV) $u(x,\cdot)$, and
for each realisation  $\omega\in\Omega$ one has an instance of a spatial field $u(\cdot,\omega)$.
To make things simple, assume that $u\in\Lp_2(\C{X}\times\Omega)$, which is isomorphic
to the tensor product $\Lp_2(\C{X})\otimes \Lp_2(\Omega) \cong 
\Lp_2(\C{X}\times\Omega)$.  Now one may investigate the splitting $p=x, \C{P}=\C{X}$
and $r(p) = u(p,\cdot)$ with $\C{U}= \Lp_2(\Omega)$ and $\C{Q}=\Lp_2(\C{X})$,
where for each $p\in\C{X}$ the model $r(p)$ is a RV.  Then the operator $C$
is on $\C{U}=\Lp_2(\Omega)$, and one usually investigates $C_{\C{Q}}$ on
$\C{Q}=\Lp_2(\C{X})$, an operator on a spatial field.  Turning everything around,
one may investigate 
the splitting $p=\omega, \C{P}=\Omega$ and $r(p)=u(\cdot,p)$ with $\C{U}=\Lp_2(\C{X})$
and $\C{Q}=\Lp_2(\Omega)$, where for each $p\in\Omega$ the model $r(p)$ is a spatial field.
The operator $C$ on $\C{U}=\Lp_2(\C{X})$ is what was before the operator $C_{\C{Q}}$
and vice versa.

\subsection{Examples and interpretations}  \label{SS:xmpl}
\begin{enumerate}
\item
 Taking up this first example, assume that the Hilbert space $\C{U}$
 is a space of centred (zero mean) random variables (RVs),
 e.g.\  $\C{U}=\Lp_2(\Omega)$ with inner product $\bkt{\xi}{\eta}_{\C{U}} := \EXP{\xi \eta}
 = \int_\Omega \xi(\omega) \eta(\omega)\,\D{P}(\di \omega)$,
 the covariance, and $r$ is a zero-mean scalar random field and $r(p)=u(p,\cdot)$
 is a zero-mean RV, or a ($n=1$) stochastic process indexed by $p\in\C{P} \subseteq \D{R}^n$.
 Then $R:\C{U}=\Lp_2(\Omega)\to\C{Q}=\Lp_2(\C{X})$ maps the RV $\xi$ to
 its spatial covariance with the random field, $R \xi = \left(p\mapsto
 \EXP{\xi(\cdot) u(p,\cdot)}\right)$.  The representation operator $R^*$ maps
 fields into random variables, $R^* v = \int_{\C{X}} v(x) u(x,\cdot)\,\di x$.
 The operator $C$ on $\C{U}=\Lp_2(\Omega)$ is
 rarely investigated, more typically one looks at $C_{\C{Q}}$ on $\C{Q}=\Lp_2(\C{X})$,
 represented by its kernel $\vkappa$ as an integral equation \refeq{eq:XVIII-int-ev}
 on the spatial domain $\C{X}$.   The kernel is the usual covariance function $\vkappa(p_1,p_2)$.
 This is the application where the name \KL{} expansion was originally used.
 We have used it here in a more general fashion.
\item Similar to the previous example, but the random field is assumed to be
 stochastically homogeneous, which means that the  covariance function
 $\vkappa(p_1,p_2)$ is \emph{shift invariant} or \emph{translation invariant}.
 This example has already been shortly described at the end of \refSS{kernel-spec-dec},
 and there is much literature on this subject, e.g.\ \citep{Karhunen1946, Karhunen1950,
 Yaglom-62-04, gelfand64-vol4, Yaglom-68-I, Yaglom-68-II, LoeveII, kreeSoize86, Matthies_encicl},
 so we will not further dwell on this.
\item 
  Here we look at the second example's interpretation of the random field described above.
  This is a special case of what has been described at the beginning of \refS{correlat},
  where the   measure $\vpi$ on $\C{P}$ is the probability measure $\vpi = \D{P}$.
  For simplicity let $r(p)$ be a centred $\C{U}$-valued RV, and each $r(p)=u(\cdot,p)$
  is an instance of a spatial field.  The associated linear operator $R:\C{U}=\Lp_2(\C{X})\to
  \C{Q}=\Lp_2(\Omega)$ maps spatial fields to \emph{weighted averages}, a RV;
  $R v = \int_{\C{X}} v(x) u(x,\cdot)\,\di x\in\C{Q}=\Lp_2(\Omega)$.  It is what
  $R^*$ was in the first example.  And here the representation operator $R^*$
  is what $R$ was in the first example.  Then $C$ is the covariance operator,
  operating on spatial fields.  This was $C_{\C{Q}}$ in the first example.
\item If $\C{P}=\{1,2,\ldots,n\}$, then $\C{U}=\spn\{r(\C{P})\}\cong \D{R}^n$
   is finite dimensional, and $\D{R}^{\C{P}} = \D{R}^n$ by definition.  Hence both 
   $\C{R} =\D{R}^n$ and $\C{Q}=\D{R}^n$, possibly with different inner products.
   In any case, $\vkappa$ is the Gram matrix of the vectors $\{r_1,\ldots,r_n\}$.
   The SVDs of $R$ are matrix SVDs, and the representation map $R^*$ is connected
   to the \KL{} expansion, which here is called the \emph{proper orthogonal decomposition} (POD).
\item If $\C{P}=[0,T]$ and $r(t) (t\in [0,T])$ is the response of a dynamical
   system with state space $\C{U}$, one may take $\C{Q}=\Lp_2([0,T])$.  The
   associated linear map $R$ tells us the dynamics of certain
   components.  To illustrate this, assume for the moment that $\C{U}=\D{R}^n$,
   a dynamical system with $n$ degrees of freedom.  Taking each canonical unit vector
   $\vek{e}_j$ in turn, one sees that $R \vek{e}_j = (t\mapsto \vek{e}_j^{\ops{T}}\vek{u}(t)
   = u_j(t))$, i.e.\ the time evolution of the $j$-th component.  The representation
   operator $R^*: \C{Q}\subset \Lp_2([0,T])$ maps scalar time-functions on their weighted
   average with the dynamics $R^* \phi = \int_{[0,T]} \phi(t) u(t)\,\di t \in \C{U}$.
\item Combining the two previous examples gives the method of \emph{temporal}
   snapshots, and the \KL{} expansion becomes the POD for a dynamical system.
\item If $\C{P}=\{ \omega_s|\; \omega_s \in \Omega\}$ are samples from some
   probability space $\Omega$, then one obtains the POD method for samples for some
   $\C{U}$-valued RV.
\end{enumerate}

\subsection{Decompositions and model reduction}  \label{SS:mod-red}
Let us go back to the example at the beginning in \refS{intro}, where
the quantity of interest is the time evolution of a dynamical system,
$t\mapsto v(t,q)$ with state space $\C{V}$, dependent on a parameter $q\in\C{S}$.
Assume for simplicity that the whole process can be thought of as an element
of $\C{V} \otimes \Lp_2([0,T])\times \C{S}) \cong \C{V} \otimes \Lp_2([0,T]) \otimes
\Lp_2(\C{S})$.  One may take $\C{U} = \C{V}\otimes\Lp_2([0,T])$, the time-histories
in state space, and $p=q$, $\C{P}=\C{S}$, and $\C{Q}=\Lp_2(\C{S})$.

But it is also possible to take $\C{U}=\C{V}$ and $p=(t,q)$, $\C{P}=[0,T]\times\C{S}$,
$\C{Q}=\Lp_2([0,T]) \otimes\Lp_2(\C{S})$.  Staying with the latter split,
for example the representation \refeq{eq:XIV} in \refT{1st-spec-rep}
becomes
\begin{equation}   \label{eq:XV}
  r(p) =  \sum_m \lambda_m^{1/2} \, s_m(p) v_m =  \sum_m \lambda_m^{1/2} \, s_m((t,q)) v_m.
\end{equation}
Now each of the scalar function $s_m((t,q))$ may be seen as a parametric model
$q\mapsto s_m(\cdot,q)$ of time functions in $\Lp_2([0,T])$.  So now one may
investigate the parametric model based on $\C{U}_* = \Lp_2([0,T])$ and
$\C{Q}_* = \Lp_2(\C{S})$ for each of the $s_m$.

Frequently the parameter space is a product space
\[
  \C{S} = \C{S}_{I} \times \C{S}_{II} \times \ldots = \prod_K \C{S}_K,\quad K=I, II, \dots ,
\]
 with a product measure $\vpi = \vpi_I \otimes \vpi_{II} \dots $, with $s_m(t,q)=
 s_m(t,(q_I,q_{II},\dots))$.  As then
\[
  \Lp_2(\C{S}) = \Lp_2(\prod_K \C{S}_K) = \Lp_2(\C{S}_I) \otimes \Lp_2(\C{S}_{II}) \otimes \dots
  = \bigotimes_K \Lp_2(\C{S}_K), \quad K=I, II, \dots ,
\]
one obtains
\begin{equation}  \label{eq:XXIII-tens}
  \C{Q} = \C{U}_*\otimes \C{Q}_* = \C{U}_* \otimes \C{Q}_I \otimes \C{Q}_{II} \otimes \dots ,
\end{equation}
with $\C{Q}_K = \Lp_2(\C{S}_K)$ for $K=I, II, \dots$.
It is clear that $\C{Q}_* = \bigotimes_K \C{Q}_K$ may be further split by different associations
depending on the value of $J$:
\begin{equation}  \label{eq:XXIII-tens-p}
  \C{Q}_* = \C{U}_{**}\otimes\C{Q}_{**} =
  \left(\bigotimes_{K=I}^J \C{Q}_K\right) \otimes \left(\bigotimes_{K>J} \C{Q}_K\right) .
\end{equation}
As will be seen, this leads to hierarchical tensor approximations, e.g.\ 
\citep{Hackbusch_tensor, boulder:2011}.

The model space has now been decomposed to
\begin{equation}  \label{eq:XXIII-tens-p-all}
  \C{U}\otimes\C{Q} = \C{U} \otimes \C{U}_* \otimes \C{Q}_* =
  \C{U} \otimes \C{U}_* \otimes
    \left(\bigotimes_K \C{Q}_K\right) .
\end{equation}

Computations usually require that one chooses finite dimensional subspaces and
bases in there, in the example case of \refeq{eq:XXIII-tens-p-all} assume that these are
\begin{align}  \label{eq:XXXa}
\text{span }\{ u_n\}_{n=1}^N =\C{U}_N \subset\C{U},\quad &\dim \C{U}_N = N,\\
  \label{eq:XXXb}
\text{span } \{\tau_k\}_{j=1}^J =\C{U}_{*J} \subset L_2([0,T]) = \C{U}_*,
  \quad &\dim \C{U}_{*J} =J,\\
  \nonumber   \forall \ell_k=1,\ldots,L_K, \; K=I, II, \dots :&\\
    \label{eq:XXXc}
\text{span } \{s_{\ell_k}\}_{\ell_k=1}^{L_k} = \C{Q}_{K, L_K}
   \subset L_2(\C{S}_K) = \C{Q}_{K}, \quad &\dim \C{Q}_{K,L_K} = L_K.
\end{align}
An approximation to $u\in \C{U}\otimes\C{Q}$ in the space described in
\refeq{eq:XXIII-tens-p-all} is thus given by
\begin{equation}  \label{eq:XXXI}
  u(x,t,q_I, \ldots, ) \approx  \sum_{n=1}^N \sum_{j=1}^J
      \sum_{\ell_I=1}^{L_I} \ldots \sum_{\ell_K = 1}^{L_K}\dots
   \tns{u}^{\ell_I,\ldots,\ell_K,\dots}_{n,j} w^n(x) \tau^j(t) 
   \left(\prod_{K} s_{\ell_K}(\omega_m)\right).
\end{equation}
Via the \refeq{eq:XXXI} one sees that the tensor
\begin{equation}  \label{eq:XXIII-tens-concr}
  \tnb{u} = \left(\tns{u}^{\ell_I,\ldots,\ell_K,\dots}_{n,j}\right)
  \in \D{R}^{(N \times J \times \prod_K L_K)}
  \cong \D{R}^N \otimes \D{R}^J \otimes \bigotimes_K \D{R}^{L_K}
\end{equation}
represents the total parametric response $u(x,t,q_I, \ldots, )$.

One way to perform model reduction is to apply the techniques described
before 
on the finite dimensional approximation space of the one described in
\refeq{eq:XXIII-tens-p-all}
\begin{equation}  \label{eq:XXIII-tens-p-all-f}
    \C{U}_N \otimes \C{Q}_M = \C{U}_N \otimes \C{U}_{*J} \otimes
    \left(\bigotimes_K \C{Q}_{K,L_K}\right) = 
    \subset  \C{U} \otimes \C{U}_* \otimes
    \left(\bigotimes_K \C{Q}_K\right) ,
\end{equation}
with $\C{Q}_{M} = \C{U}_{*J} \otimes \left(\bigotimes_K \C{Q}_{K,L_K}\right)$
and dimension $\dim \C{Q}_M = M= J \times \prod_K L_K$
but not using the full dimension, as the spectral analysis of the `correlation'
operator $C$ picks out the important parts.

Another kind of reduction works directly with the tensor $\tnb{u}$
in \refeq{eq:XXIII-tens-concr}.  It has formally
$R^\prime= N \times J \times \prod_K L_K$ terms.
Here we only touch on this subject, which is a \emph{nonlinear} kind
of model reduction, and that is to represent this tensor with many
times fewer $R \ll R^\prime$ terms through
what is termed a low-rank approximation, for a thorough treatement
see the monograph \citep{Hackbusch_tensor}.

Whereas the so called \emph{canonical polyadic} (CP) decomposition 
uses the flat tensor product in \refeq{eq:XXIII-tens-p-all}
--- under the name \emph{proper generalised decomposition} (PGD)
\citep{Nouy2009, NouyACM:2010, ammChin2010, LadevezeCh2011} this is
also a computational method to solve an effectively high-dimensional
problem as  \refeq{eq:I} or \refeq{eq:I-p}, see the review
\citep{chinestaPL2011} and the monograph \citep{chinestaBook} ---
the recursive use of splittings \refeq{eq:XXIII-tens-p} leads to
\emph{hierarchical} tensor approximations, e.g.\ \citep{Grasedyck2010}.
  The index set can be
thought to be partitioned and arranged into a binary tree, each time
causing a split as in \refeq{eq:XXIII-tens-p}, or rather on the finite
dimensional approximation \refeq{eq:XXIII-tens-p-all-f},
or equivalently in the concrete tensor in \refeq{eq:XXIII-tens-concr}.
Particular cases of this are the \emph{tensor train} (TT)
\citep{oseledetsTyrt2010, oseledets2011} and more
generally the \emph{hierarchical Tucker} (HT) decompositions,
see the review \citep{GrasedyckKressnerTobler2013} and the
monograph \citep{Hackbusch_tensor}.  An example how this representation
then allows also fast post-processing such as finding maxima is given
in \citep{Espig:2013}.
Let us also mention that these sparse or low-rank tensor representations
are connected with the expressive power of deep neural networks
\citep{CohenSha2016, KhrulkovEtal2018}.  Neural networks are one
possibility of choosing the approximation functions in
\refeq{eq:XXXc}.
Obviously, a good reduced order model is one with only few terms.
One recognises immediately that the SVD structure of the associated
linear map of such a split determines how many terms are needed for
a good approximation.  Equivalently it is the structure of the spectrum of the
appropriate correlation operator associated with the splitting.

%
%
%
%
%
%
%
%
%
%
%
%
%


%

\section{Refinements}  \label{S:refine}
Often the parametric element has more structure than is resolved by saying that
for each $p\in\C{P}$  one has $r(p)$ in some Hilbert space $\C{U}$.  Most of
the preceding had to do with alternative ways of describing the dependence on the
parameter $p$.  Here a short look is taken on the case when the Hilbert space $\C{U}$
has more structure, which one might want to treat separately.  One big area,
which we only entered slightly, are invariance properties as the invariance
w.r.t.\ shifts for stationary or stochastically homogeneous random fields touched
on in \refSS{kernel-fact}.  We shall look only at two simple but
instructive cases.

\subsection{Vector fields}   \label{SS:vector_fields}
One of the simplest variations on the modelling in the previous sections is
the refinement that the  r\^ole of the Hilbert space $\C{Q}$ is taken by a
tensor product $\C{W}=\C{Q} \otimes \C{E}$, where as before $\C{Q}$ is a 
Hilbert space of scalar functions and  $\C{E}$ a \emph{finite-dimensional}
inner-product (Hilbert) space \citep{kreeSoize86}.
The associated linear map is then a map
\begin{equation}   \label{eq:vect-R}
   R_{\C{E}}: \C{V} \to \C{W} = \C{Q}\otimes \C{E}.
\end{equation}
One possible situation where this occurs is when, similar to the
third example in \refSS{xmpl}, the random field $\vek{u}(x,p)$ is not scalar-
but vector valued, i.e.\ $\vek{u}(x,p) \in \C{E}$.  It could be that several
correlated scalar fields have to be described which have been collected into
a vector $[u_a(x,p), \dots, u_j(x,p)]$, or that $\vek{u}(x,p) \in \D{R}^n$ is
actually to be interpreted as a vector in the space $\D{R}^n$, e.g.\ a velocity
vector field.  Without loss of generality we shall assume that $\C{E} = \D{R}^n$.
This obviously also covers the case when $\C{E}$ is a space of tensors of higher
degree; although for tensors of even degree we shall show a further simplification
in \refSS{tensor_fields}.

In this case, when $\C{V} = \C{U} \otimes \C{E}$, the parametric map is
\begin{equation}   \label{eq:param-vect-r}
   \tns{r}:\C{P}\to\C{V} = \C{U} \otimes \C{E};\quad \tns{r}(p) = \sum_k r_k(p) \vek{r}_k,
\end{equation}
where as before $r_k(p) \in \C{U}$ --- here in the motivating example a Hilbert
space of scalar fields ---
and the $\vek{r}_k \in \C{E}$.   In this case the associated map $R_{\C{E}}$
is chosen \citep{kreeSoize86} to be 
\begin{equation}   \label{eq:vect-R-ex}
   R_{\C{E}} = \sum_k R_k\otimes\vek{r}_k:\C{U} \ni u \mapsto 
   \sum_k R_k(u) \vek{r}_k = \sum_k \bkt{u}{r_k(\cdot)}_{\C{U}}\vek{r}_k \in \C{Q}\otimes \C{E}.
\end{equation}
where the maps $R_k:\C{U}\to\C{Q}$ are just the maps from \refeq{eq:IV},
but each $R_k$ is the associated map to $r_k(p)$.

The `correlation' can now be given by a bilinear form, but not with
values in $\D{R}$ as in \refD{corr}, but with
values in $\C{E}\otimes\C{E}$.  For this we define on $\C{W}^2
= (\C{Q}\otimes\C{E})^2$ a bilinear form $[\cdot\mid\cdot]$ with values
in $\C{E}\otimes\C{E}$ first on elementary tensors,
\begin{equation}   \label{eq:vect-biform}
\forall \tns{s} = s\otimes\vek{s}, \tns{t}= \tau\otimes\vek{\tau} \in \C{W}=\C{Q}\otimes\C{E}:\quad
   [s\otimes\vek{s}\mid \tau\otimes\vek{\tau}] := \bkt{s}{\tau}_{\C{Q}}\, \vek{s}\otimes\vek{\tau} , 
\end{equation}
and then extended by linearity.
Concerning $\C{U}$ and $\C{Q}$ we make the
same assumptions as before in \refSS{ass-lin-map} and \refS{correlat}.

\begin{defi}[Vector-Correlation]  \label{D:corr-vec}
Define
 a densely defined map $C_{\C{E}}$ in $\C{V}=\C{U}\otimes\C{E}$  on elementary tensors as
\begin{multline}   \label{eq:vect-corr}
\forall \tns{u} = u\otimes\vek{u}, \tns{v}=v\otimes\vek{v} \in \C{V}=\C{U}\otimes\C{E}:\\
   \bkt{C_{\C{E}}\tns{u}}{\tns{v}}_{\C{U}} := \vek{u}^{\ops{T}}\,[ R_{\C{E}} u \mid R_{\C{E}}v]\,
   \vek{v} =    \sum_{k,j}  \bkt{R_k(u)}{R_j(v)}_{\C{Q}}
    \,(\vek{u}^{\ops{T}}\vek{r}_k)\, (\vek{r}_j^{\ops{T}}\vek{v})
\end{multline}
and extend it by linearity.  It may be called the 
\emph{`correlation'} operator.  By construction it is self-adjoint
and positive.
\end{defi}

The factorisations and decompositions then have to be of this operator.
The eigenproblem on $\C{V}$ is: Find $\lambda\in\D{R}$, $\tns{v}\in\C{U}\otimes\C{E}$
with $\tns{v}= \sum_\ell v_\ell \otimes\vek{r}_\ell$, such that
\begin{equation}   \label{eq:vect-eig-U}
   C_{\C{E}}\tns{v}  = \sum_{k,j,\ell} (R_j^* R_k v_\ell)\, (\vek{r}_k^{\ops{T}}\vek{r}_\ell)\,
   \vek{r}_j  = \lambda \tns{v} = \lambda  \sum_j v_j \otimes\vek{r}_j .
\end{equation}

The kernel $\vek{\vkappa}_{\C{E}}: \C{P}^2 \to (\C{E}\otimes\C{E})$ for the
eigenvalue problem on $\C{W}=\C{Q}\otimes\C{E}$ is
\begin{equation}   \label{eq:vect-kernel}
   \vek{\vkappa}_{\C{E}}(p_1,p_2) = \sum_{k,j} \bkt{r_k(p_1)}{r_j(p_2)}_{\C{V}} 
   \; \vek{r}_k\otimes\vek{r}_j  .
\end{equation}
So $\vek{\vkappa}$ is matrix-valued.  These are
actually `correlation' matrices.

In case the function space $\C{Q}$ has the structure of $\Lp_2(\C{P})$ with
measure $\vpi$ on $\C{P}$, the Fredholm eigenproblem has the following form:
Find $\lambda\in\D{R}$, $\tns{s}\in\C{Q}\otimes\C{E}$ with $\tns{s}=
\sum_\ell \vsigma_\ell(\cdot)\vek{r}_\ell$ such that
\begin{multline}   \label{eq:vect-Fredh-eig}
   \int_{\C{P}} \vek{\vkappa}_{\C{E}}(p_1,p_2) \left(\sum_\ell \vsigma_\ell(p_2)
   \vek{r}_\ell\right) \,  \vpi(\di p_2)=\\
   \sum_{k,j,\ell} \left( \int_{\C{P}} \bkt{r_k(p_1)}{r_j(p_2)}_{\C{U}}\, \vsigma_\ell(p_2) \;
   \vpi(\di p_2)    \right)(\vek{r}_k^{\ops{T}}\vek{r}_\ell) \vek{r}_j  = \lambda \sum_j 
   \vsigma_j(p_1) \vek{r}_j .
\end{multline}

Both of these eigenproblems then combine into a generalised \KL{} expansion,
the analogue of \refeq{eq:XIV} in \refT{1st-spec-rep}:
\begin{equation}   \label{eq:vect-KL}
   \tns{r}(p) = \sum_k r_k(p) \vek{r}_k = \sum_m \lambda_m^{1/2} 
   \left( \sum_k \vsigma_{m,k}(p) v_{m,k} \,\vek{r}_k \right) = \sum_k \left( 
   \sum_m \lambda_m^{1/2} \vsigma_{m,k}(p) v_{m,k} \right) \,\vek{r}_k.
\end{equation}

\subsection{Tensor fields}   \label{SS:tensor_fields}
Some situations as described in the previous \refSS{vector_fields} allow
an even somewhat simpler approach.  This is the case when the vector space
$\C{E}$ in \refeq{eq:vect-R} consist of tensors of \emph{even} degree.
Formally this means that $\C{E} = \C{F}\otimes\C{F}$ for some space
of tensors $\C{F}$ of half the degree.
A tensor of even degree can always be thought of as a linear map from a
space of tensor of half that degree into itself.  Namely, let for example
$\tensor*{\tnb{A}}{_a_b_c^d^e^f}\in \C{E}$ be a tensor of \emph{even} degree---here six.
Then this tensor acts as a linear map on the space of tensors of e.g.\ the
form $\tensor*{\tnb{f}}{^b_e_f}\in\C{F}$ (the Einstein summation convention for
tensor contraction is used in this symbolic index notation):
\[
   \tensor*{\tnb{A}}{_a_b_c^d^e^f}   \tensor*{\tnb{f}}{^b_e_f} = \tensor*{\tnb{q}}{^d_a_c}.
\]
Often the particular application domain will dictate which space of tensors it
acts on.
Being a linear map, it can be represented as a \emph{matrix} $\vek{A}\in\D{R}^{n\times n}$,
which we shall assume
from now on.  Often these linear maps / matrices have to satisfy some additional properties,
for example they have to be positive definite or orthogonal.

It is maybe now the opportunity to make an important remark:
The representation methods which have been shown here are \emph{linear} methods,
which means they work best when the object to be represented is in a 
\emph{linear} or \emph{vector space}, essentially free from nonlinear constraints.
Consider two illustrative examples:

As a first example, assume that $\vek{A}$ has to be orthogonal, then one requires
$\vek{A}^{\ops{T}}\vek{A} =\vek{I}=\vek{A}\vek{A}^{\ops{T}}$, a nonlinear constraint.
It is well known that the
orthogonal matrices $\ops{O}(n)$ form a compact Lie group, just as the sub-group
of special orthogonal matrices $\ops{SO}(n)$.  Here it is important to notice that
their Lie algebra $\F{o}(n) = \F{so}(n)$, the tangent space at the group identity $\vek{I}$,
are the \emph{skew} symmetric matrices, a \emph{free} linear space.
And each $\vek{Q}\in \ops{SO}(n)$ can be represented with the exponential map
$\vek{Q} = \exp(\vek{S})$ with $\vek{S}\in\F{so}(n)$ and `$\exp(\cdot)$' the matrix exponential.
This recipe, using the exponential map from the Lie algebra,
which is a vector space, to its corresponding
Lie group, is a very general one.  One has to deal only with representations on free linear
spaces, the Lie algebra, but models entities in the Lie group.

For another example, assume that the matrix $\vek{A}\in\ops{Sym}^+(n)$ has to be symmetric positive
definite (spd), as is often required when one wants to model constitutive material tensors.
One defining condition is that it can be
factored as $\vek{A}=\vek{G}^{\ops{T}}\vek{G}$ with invertible $\vek{G}\in\ops{GL}(n)$.
Both of these are nonlinear constraints.  In fact the spd matrices
are an open cone, a Riemannian manifold, in the space of
all symmetric matrices $\F{sym}(n)$.
There are different ways how to make $\ops{Sym}^+(n)$ into a Lie group, but the
important thing here is that any $A\in\ops{Sym}^+(n)$ can be represented
again with the matrix exponential as
$\vek{A}=\exp(\vek{H})$ with $\vek{H}\in\F{sym}(n)$, a \emph{free} linear space.
Let us point out that this also important in the case $n=1$, i.e.\ when $\vek{A}$ is
a positive scalar.  Here $\F{sym}(1)=\D{R}$ and the map $\exp(\cdot)$ is the usual exponential.

A parametric element in this special case of \refeq{eq:param-vect-r}, let us say in
the example of positive definite matrices $\vek{A}(p)\in\C{Q}\otimes\ops{Sym}^+(n)$,
would be represented by an element $\vek{H}(p)\in\C{Q}\otimes\F{sym}(n)$
and then exponentiated:
\begin{equation}   \label{eq:exp-rep-spd}
   \vek{H}(p) = \sum_k \vsigma_k(p) \vek{H}_k,\qquad
   \vek{H}(p) \mapsto \exp(\vek{H}(p)) = \vek{A}(p) .
\end{equation}
This way one is sure that $\vek{A}(p)\in\ops{Sym}^+(n)$ for each $p\in\C{P}$.
Therefore we can now concentrate on the problem of representing $\vek{H}(p)$.

Here everything is very similar to the previous \refSS{vector_fields}.
The associated linear map in \refeq{eq:vect-R} remains as it is, only
that now $\C{E}=\F{sym}(n)$.  The parametric map would be written as
\begin{equation}   \label{eq:tens-r-param}
   \tns{R}(p) = \sum_k r_k(p)\otimes\vek{R}_k \in \C{U}\otimes \C{E},\quad
   \text{ with } \quad \vek{R}_k \in \F{sym}(n).
\end{equation}
The correlation corresponding to
\refD{corr-vec} is now defined as

\begin{defi}[Tensor-Correlation]  \label{D:corr-tens}
Define
 a densely defined map $C_{\C{E}}$ in $\C{W}=\C{U}\otimes\C{F} = \C{U}\otimes\D{R}^n$ 
 --- observe, \emph{not} $\C{U}\otimes\C{E}=\C{U}\otimes\F{sym}(n) = 
 \C{U}\otimes\D{R}^{n(n+1)/2}$ --- on elementary
 tensors through a bilinear form  as
\begin{multline}   \label{eq:tens-corr}
\forall (\tns{u} = u\otimes\vek{v}), (\tns{v}=v\otimes\vek{v}) \in \C{W}=\C{U}\otimes\C{F}:\\
   \bkt{C_{\C{F}}\tns{u}}{\tns{v}}_{\C{U}} :=  \sum_{k,j}  \bkt{R_k(u)}{R_j(v)}_{\C{Q}}
    \,(\vek{R}_k\vek{u})^{\ops{T}}(\vek{R}_j \vek{v}),
\end{multline}
and extend it by linearity.  It may be called the 
\emph{`correlation'} operator.  By construction it is self-adjoint
and positive.
\end{defi}

The eigenproblem for $C_{\C{F}}$ corresponding to \refeq{eq:vect-eig-U} is now formulated on
$\C{W}=\C{U}\otimes\C{F} = \C{U}\otimes\D{R}^n$, i.e.\ like \refeq{eq:vect-eig-U}
but with $\C{E}$ replaced by $\C{F}$, otherwise everything is as before
and hence will not be spelled out in detail.
The eigenvectors, analogous to \refeq{eq:vect-eig-U}, will look the same as there but with
$\vek{v}_{m,j}\in\C{F}$:
\begin{equation}   \label{eq:eig-tens-U}
   \C{W}\ni \tns{v}_m = \sum_j v_{m,j} \otimes \vek{v}_{m,j} \in \C{U}\otimes\C{F} .
\end{equation}

The kernel corresponding to \refeq{eq:vect-kernel} is
\begin{equation}   \label{eq:tens-kernel}
   \vek{\vkappa}_{\C{F}}(p_1,p_2) = \sum_{k,j} \bkt{r_k(p_1)}{r_j(p_2)}_{\C{U}} 
   \; \vek{R}_k^{\ops{T}} \vek{R}_j  ,
\end{equation}
with eigenvectors of the form $\tns{s}_m(p) = \sum_j \vsigma_{m,j}(p) \vek{v}_{m,j}
\in\C{Q}\otimes\C{F}$.  So $\vek{\vkappa}_{\C{F}}$ is matrix-valued here as well.
From these the final representation, 
the analogue of \refeq{eq:XIV} in \refT{1st-spec-rep} and \refeq{eq:vect-KL}, is obtained as:
\begin{equation}   \label{eq:tens-KL}
   \tns{R}(p) = \sum_k r_k(p) \vek{R}_k = \sum_m \lambda_m^{1/2} 
   \left( \sum_j \vsigma_{m,j}(p) v_{m,j} \,\vek{v}_{m,j} \otimes \vek{v}_{m,j} \right) .
\end{equation}

%
%
%
%
%
%
%
%
%
%
%
%
%


%

\section{Conclusion} \label{S:concl}
A parametric mapping $r:\C{P}\to\C{U}$ has been analysed in a  variety of
settings.  The basic idea is the associated linear map $R:\C{U}\to\D{R}^{\C{P}}$,
which both generalises the parametric mapping and enables the linear analysis.
This leads immediately to the RKHS setting, and a first equivalent
representation on the RKHS in terms of tensor products.  By choosing other inner
products than the one coming from the RKHS one can analyse the importance of different
features, i.e.\ subsets of the parameter set.

Importance is measured by the spectrum of the `correlation' operator $C$,
and the representation is again in terms of tensor products.  The correlation
is factored by the associated linear map, and it is shown that on one hand
all factorisations are unitarily equivalent, and on the other hand that each
factorisation leads to
differently parametrised representations, indeed linear resp.\ affine representations,
if the correlation is a nuclear or trace class operator.
In fact, each such factorisation corresponds to a representation, and vice versa.
This equivalence is due to our strict assumptions on the associated linear map,
but in general the associated linear map is a truly more general concept.

These linear representations are in terms of real valued functions, which may be seen
as some kind of `co-ordinates' on the parameter set.  In the RKHS case, they are truly
co-ordinates or an embedding of the parameter set into the RKHS, as each
parameter point $p\in\C{P}$ can be identified with the evaluation functional
$\updelta_p$ and hence the kernel $\vkappa(p,\cdot)$ at that point.

An equivalent spectral analysis can be carried out for the kernel space in
terms of integral equations and integral transforms.  These spectral decompositions
or other factorisations also lend themselves to the construction of parametric
reduced order models, as the importance of different terms in the representing
\KL- or POD-series can be measured.  But other factorisations also lead to,
not necessarily optimal, reduced order models.  For the sake of simplicity
for the representation only orthonormal bases have been considered here as they
appear quite natural in the Hilbert space setting, but obviously other
bases can be considered.  The tensor product nature
of this series makes it natural to employ this factorisation in a recursive
fashion and thereby to generate a representation through high-order tensors.
These often allow very efficient low-rank approximations, which is in fact another,
but this time nonlinear, model order reduction.
Certain refinements are possible in case the representation space has the structure
of a `vector'- or `tensor'-field, a point which is only briefly touched.
It was also shown that the structure of the spectrum of the correlation operator
attached to such a tensor-space factorisation, or equivalently the structure
of the set of singular values of the associated linear map, determines how many
terms are needed in a reduced order model to achieve a certain degree of accuracy.

Thus the functional analytic view on parametric problems via 
decompositions of linear maps
 gives a certain unity to seemingly different procedures
which turn out to be closely related, at least if one looks for the
similarities.  This constitutes a natural introduction and background to low-rank
tensor product representations, which are crucial for efficient
computation.  They are naturally employed in a \emph{functional approximation}
approach to parametric problems.

\ignore{
\subsection{Test of math fonts} \label{SS:testf}
\[
\E{D}(\Omega) \hookrightarrow \E{S}(\Omega) \hookrightarrow \E{E}(\Omega)  \hookrightarrow
   \Lp_1^{\text{loc}}(\Omega) \hookrightarrow \E{E}^\prime(\Omega) \hookrightarrow
   \E{S}^\prime(\Omega) \hookrightarrow \E{D}^\prime(\Omega) \hookleftarrow
\]

\paragraph{Slanted:} --- serifs; \emph{for variables, scalars and vectors, operators}
\[ a, \alpha, A, \phi, \vphi, \pi, \theta, \vtheta, \delta, \Phi, \Pi, \Theta, \Delta; \qquad 
   \vek{a}, \vek{\alpha}, \vek{A}, \vek{\phi}, \vek{\vphi}, \vek{\pi}, \vek{\theta},
   \vek{\vtheta}, \vek{\delta}, \vek{\Phi}, \vek{\Pi}, \vek{\Theta}, \vek{\Delta}. \]

\paragraph{Slanted:} --- sans serif; \emph{for variables, tensors of higher degree}
\[ \tns{a}, \tns{\alpha}, \tns{A}, \tns{\phi}, \tns{\vphi}, \tns{\pi}, \tns{\theta},
   \tns{\vartheta}, \tns{\delta}, \tns{\Phi}, \tns{\Pi}, \tns{\Theta}, \tns{\Delta}; \qquad
   \tnb{a}, \tnb{\alpha}, \tnb{A}, \tnb{\phi}, \tnb{\vphi}, \tnb{\pi}, \tnb{\theta},
   \tnb{\vartheta}, \tnb{\delta}, \tnb{\Phi}, \tnb{\Pi}, \tnb{\Theta}, \tnb{\Delta}. \]
   
\paragraph{Other:}  --- no small letters and no Greek letters;
      \emph{for sets, spaces, special operators}
\[ \C{E}, \C{Q}, \C{R};\qquad \D{E}, \D{Q}, \D{R};\qquad \E{E}, \E{Q}, \E{R}.\]

\paragraph{Fraktur:} --- no Greek letters; \emph{for logical variables, multi-indices?, sets of sets}
\[ \F{e}, \F{q}, \F{r};\qquad \F{E}, \F{Q}, \F{R}. \]

\paragraph{Upright:} --- serifs and sans-serif; \emph{for constants, special operators}
\[ \mrm{a}, \mrm{A}, \mrm{\Phi}, \mrm{\Pi}, \mrm{\Theta}, \mrm{\Delta};\qquad
   \mat{a}, \mat{A}, \mat{\Phi}, \mat{\Pi}, \mat{\Theta}, \mat{\Delta};\qquad 
   \ops{a}, \ops{A}, \ops{\Phi}, \ops{\Pi}, \ops{\Theta}, \ops{\Delta};\qquad
   \opb{a}, \opb{A}, \opb{\Phi}, \opb{\Pi}, \opb{\Theta}, \opb{\Delta}; \]
 --- upright Greek letters, especially small ones, and alternative capital ones
\[ \upalpha, \upphi, \upvarphi, \uppi, \uptheta, \upvartheta, \updelta,
   \Upphi, \Uppi, \Uptheta, \Updelta; \qquad
   \vek{\upalpha}, \vek{\upphi}, \vek{\upvarphi}, \vek{\uppi}, \vek{\uptheta}, \vek{\upvartheta},
   \vek{\updelta}, \vek{\Upphi}, \vek{\Uppi}, \vek{\Uptheta}, \vek{\Updelta}.  \]
}

%
%
%
%
%
%
%
%
%
%
%
%
%





\bibliography{\thebib/jabbrevlong,\thebib/stochastic,\thebib/fuq-new,%
\thebib/fa,\thebib/mat_BU-1-S,\thebib/phys_D,\thebib/num,\thebib/highdim}


{ 
   \tiny
       \texttt{\RCSId} 
   }



\end{document}